\documentclass[twoside,11pt]{article}

\usepackage{blindtext}

\usepackage[preprint]{jmlr2e} 
\usepackage{bm}
\usepackage[title]{appendix}
\usepackage{listings}
\usepackage{amsfonts,mathtools,mathrsfs,amsmath}
\usepackage{epstopdf}
\usepackage{enumerate}
\usepackage{enumitem}
\usepackage{float}
\usepackage{xcolor}
\usepackage{color}
\usepackage{slantsc}
\usepackage{listings}
\usepackage{hyperref}
\usepackage{datetime2}
\usepackage{extarrows} 
\usepackage{booktabs, siunitx}
\usepackage{algorithm, algpseudocode}

\DeclareMathOperator*{\argmin}{argmin}
\DeclareMathOperator*{\E}{\mathbb{E}}

\DeclareMathOperator*{\Cov}{\text{Cov}}

\newcommand\R{\,{\mathbb{R}}}   
\newcommand{\T}[1]{\text{#1}}   
\newcommand{\Bra}[1]{\left({#1}\right)} 
\newcommand{\SqBra}[1]{\left[{#1}\right]} 
\newcommand{\Set}[1]{\left\{{#1}\right\}} 
\newcommand{\inner}[2]{\langle #1, #2 \rangle} 

\usepackage{lastpage}
\jmlrheading{}{}{1-\pageref{LastPage}}{}{}{}{Haotian Chen, Yu-Wei Chen, Qifan Song and Jordan Awan} 

\ShortHeadings{Privacy and Smoothing in Functional Mean Estimation}{Chen, Chen, Song and Awan}
\firstpageno{1}

\begin{document}

\title{On the Order of Privacy and Smoothing\\ in Functional Mean Estimation}

\author{\name Haotian Chen \email chen5204@purdue.edu \\
       \addr Department of Statistics\\
       Purdue University\\
       West Lafayette, IN 47907, USA
       \AND
       \name Yu-Wei Chen \email chen4357@purdue.edu \\
       \addr Department of Statistics\\
       Purdue University\\
       West Lafayette, IN 47907, USA
       \AND
       \name Qifan Song \email qfsong@purdue.edu \\
       \addr Department of Statistics\\
       Purdue University\\
       West Lafayette, IN 47907, USA
       \AND
       \name Jordan Awan \email jaa557@pitt.edu \\
       \addr Department of Statistics\\
       University of Pittsburgh\\
       Pittsburgh, PA 15260, USA}
       
\editor{}

\maketitle

\begin{abstract}
Differential privacy (DP) is a leading framework for privacy protection, yet privacy mechanisms tailored to functional data analysis (FDA) remain relatively underexplored. For functional mean estimation, we investigate a fundamental question regarding the ordering of privatization and estimation by comparing two canonical pipelines: \emph{Fit-Then-Protect} (FTP), which smooths trajectories before privatization, and \emph{Protect-Then-Fit} (PTF), which privatizes observations at discrete time points before smoothing. When data are observed without error, we prove that FTP achieves strictly smaller mean squared error than PTF asymptotically (as the sample size increases) under the same $\mu$-Gaussian differential privacy ($\mu$-GDP) budget in both finite-dimensional Hilbert spaces and infinite-dimensional reproducing kernel Hilbert spaces (RKHS). Simulation studies support our theory and explore settings with observation error and fitting constraints, demonstrating that FTP outperforms PTF under small observation noise, whereas PTF can become preferable in high-noise settings. In addition to this comparison, we establish a novel GDP guarantee for Gaussian perturbations taking values in a Banach space, extending existing results beyond the RKHS setting.
\end{abstract}


\section{Introduction}
Differential privacy (DP) provides a rigorous framework for protecting individual-level information while balancing data usability. By carefully calibrating the noise added to the released summaries, DP obscures the contribution of any particular individual. Since \citet{dwork2006calibrating} formalized DP, it has become a central paradigm in privacy-preserving data analysis, stimulating extensive research in statistics, machine learning, and related fields. Despite this progress, much of the classical DP literature has focused on database queries, finite-dimensional statistical summaries, and conventional machine-learning tasks \citep{wasserman2010statistical, chaudhuri2011differentially, dwork2014algorithmic}. Extending DP to more complex data objects raises additional theoretical and methodological challenges, particularly when the data exhibit infinite-dimensional structure and complex dependence over continuous domains.

These challenges are particularly relevant in functional data analysis (FDA), where the basic unit of interest is a function or curve, potentially infinite-dimensional, but typically observed only at finitely many discrete points. Classical examples include growth curves, wearable sensor data, and other densely or continuously recorded signals \citep{wang2016functional,ramsay2005functional,acar2025functional}. Instead of treating observations at different time points independently, FDA models each subject's data as a smooth function, allowing statistical procedures to exploit the within-curve correlation and geometric structure inherent in functional observations. Over the past several decades, FDA has become an important area of modern statistics with applications across biology, medicine, engineering, and environmental science \citep{ferraty2006nonparametric, kokoszka2017introduction}.

As functional data become increasingly prevalent, protecting the privacy of function-valued information has become an important challenge. In many applications, the object of interest is itself a function—for example, a mean trajectory, a regression function, or an estimated individual curve—which may reveal sensitive information about individuals. While classical DP mechanisms are primarily designed for finite-dimensional statistics, extending privacy guarantees to infinite-dimensional function spaces is substantially more challenging, requiring both suitable notions of sensitivity and carefully constructed noise processes. These challenges naturally raise the question of how DP should be formulated and implemented for function-valued outputs. Previous work has built $(\epsilon,\delta)$-DP for functionals in different function spaces \citep{hall2013differential, mirshani2019formal}, while recent work extends this to the Gaussian differential privacy (GDP) framework under reproducing kernel Hilbert spaces (RKHS) \citep{raisa2024noise}.

Another fundamental question arises regarding when and where privacy should be enforced, given privacy mechanisms for functional outputs are available. For a functional estimation problem, privacy can be incorporated in at least two natural ways. One may first smooth the functional trajectories and then privatize the resulting estimates, or alternatively privatize the data first and subsequently perform the estimation using the protected observations. We refer to these two strategies as \emph{Fit-Then-Protect} (FTP) and \emph{Protect-Then-Fit} (PTF), respectively. Although both procedures can be calibrated to have the same privacy guarantee, they exploit the functional structure of the data in fundamentally different ways, and little is known about when one ordering is superior. 

\textbf{Contributions:} In this paper, we study how the ordering of privacy protection and functional estimation affects the statistical accuracy in functional mean estimation. We first establish a sufficient condition under which Gaussian perturbations of function-valued summaries in general Banach spaces satisfy a prescribed GDP guarantee, extending existing results beyond the RKHS setting. We then compare the FTP and PTF pipelines under a common GDP budget. When functional observations are measured without error, we prove that FTP achieves strictly smaller estimation error than PTF asymptotically as the sample size increases, in both finite-dimensional Hilbert spaces equipped with the norm induced by $L^2[0,1]$ inner product and reproducing kernel Hilbert spaces (RKHSs). Numerical experiments support these theoretical findings and further investigate the impact of observation noise, showing that FTP performs better when observation noise is small, whereas PTF may be preferable when observation noise is large. We also examine how estimation constraints and the GDP privacy parameter $\mu$ affect the accuracy of the private estimators.

\textbf{Related Work:} Previous work on differential privacy for function-valued objects has primarily focused on extending classical privacy mechanisms from finite-dimensional vectors to infinite-dimensional function spaces. \citet{hall2013differential} proposed one of the earliest Gaussian process perturbation mechanisms for releasing functional summaries while satisfying $(\epsilon,\delta)$-DP in RKHS. Building upon this framework, \citet{mirshani2019formal} established a general theory of Gaussian perturbation mechanisms in Banach spaces and characterized sufficient conditions under which function-valued summaries satisfy $(\epsilon,\delta)$-DP. Complementary work studied the limitations of elliptical perturbations in infinite-dimensional settings \citep{reimherr2019elliptical} and extended the exponential mechanism to functional summaries, including functional principal component analysis \citep{awan2019benefits}. More recently, \citet{raisa2024noise} established GDP guarantees for a Gaussian-process-based functional mechanism in RKHS. Beyond general privacy mechanisms, \citet{xue2024optimal} studied private distributed FDA and established minimax results for functional mean and varying-coefficient estimation. 

Another line of research investigates the role of dependence and structural information in privacy mechanisms. \citet{cao2018quantifying} analyzed privacy leakage under temporal correlations in continuous data release, while other work has developed correlated-noise mechanisms to preserve time-series structure and mitigate vulnerabilities of independent perturbations \citep{wang2017cts,mcelroy2023flip}. More broadly, geometry-aware and structured mechanisms have been shown to achieve near-optimal or improved utility over conventional independent perturbations in finite-dimensional settings, including mechanisms for vector- and matrix-valued queries, private learning \citep{nikolov2013geometry, chanyaswad2018mvg, awan2021structure, choquette2024correlated}. 

\textbf{Organization:} The remainder of this paper is organized as follows. Section~\ref{sec: Background} reviews the necessary background on functional data analysis and differential privacy. Section~\ref{sec: Problem Setup} formulates the problem of differentially private functional mean estimation and introduces the \emph{Fit-Then-Protect} (FTP) and \emph{Protect-Then-Fit} (PTF) strategies. Section~\ref{sec: Theoretical Results} establishes a Gaussian differential privacy mechanism for function-valued summaries and develops theoretical comparisons between the two strategies. Section~\ref{sec: Simulation Study} presents simulation studies evaluating their empirical performance under different function spaces in various settings. Section~\ref{sec: Discussion} concludes the paper with a discussion of potential extensions. Proofs, technical details, as well as supplementary simulation results are in the appendix. The source code used to implement the simulation study and data analysis is available at \url{https://github.com/Hao-tian-Chen/DP-FDA}.

\section{Background}
\label{sec: Background}
In this section, we review the necessary background for FDA and DP. In Section~\ref{subsec: Functional Data Analysis}, we give the notation and definition for some common function spaces in FDA, such as Banach spaces and reproducing kernel Hilbert spaces (RKHS), define the Gaussian random element in a separable Banach space, and introduce the corresponding Cameron-Martin space and discuss its relationship with RKHS. In Section~\ref{subsec: Differential Privacy}, we first define $\mu$-GDP, then introduce the definition of sensitivity and the Gaussian mechanisms to ensure $\mu$-GDP, and elaborate on the GDP properties for composition and post-processing. 

\subsection{Functional Data Analysis}
\label{subsec: Functional Data Analysis}
Let $(\mathbb{B}, \|\cdot\|_{\mathbb{B}})$ be a real separable Banach space equipped with its Borel $\sigma$-algebra. A linear mapping $l: \mathbb{B} \to \mathbb{R}$ is a \emph{continuous linear functional} if it is bounded, that is, $\sup_{\|h\|_{\mathbb{B}} \le 1} |l(h)| < \infty$. The collection of all such continuous linear functionals forms the topological dual space of $\mathbb{B}$, denoted by $\mathbb{B}^*$. Equipped with the dual norm defined by $\|l\|_{\mathbb{B}^*} = \sup_{\|h\|_{\mathbb{B}} \le 1} |l(h)|$, the dual space $\mathbb{B}^*$ is itself a Banach space \citep[Proposition 1.29]{fabian2011banach}.

A $\mathbb{B}$-valued random element $Z$ is called Gaussian if $l(Z)$ is a normally distributed real random variable for every $l\in\mathbb{B}^*$ \citep{ledoux1991probability}. We write $Z\sim\mathcal{N}_{\mathbb{B}}(m,C)$, where $m\in\mathbb B$ is its mean, characterized by $\mathbb{E}[l(Z)] = l(m),\ l\in \mathbb{B}^*$,  and $C:\mathbb{B}^*\to\mathbb{B}$ is its covariance operator, defined by
\[
C l 
=
\mathbb{E}\big[l(Z-m)(Z-m)\big],
\qquad l\in\mathbb{B}^*.
\]

A Gaussian random element $Z$ is called centered if $m = 0$. In this case, we write $Z\sim \mathcal{N}_{\mathbb{B}}(0,C)$, and its covariance operator is given by $Cl = \mathbb{E}[l(Z)Z],\ l\in\mathbb{B}^*$. Unless otherwise stated, we consider centered Gaussian random elements throughout the remainder of this work.

\paragraph{Cameron--Martin Space.} Let $Z\sim \mathcal{N}_{\mathbb{B}}(0,C)$, and equip the range $C(\mathbb{B}^*) \subset \mathbb{B}$ with the inner product
\[
\inner{Cl}{Cg} := l(Cg) = \mathbb{E} [l(Z)g(Z)], \qquad l, g \in \mathbb{B}^*.
\]
The \emph{Cameron--Martin (CM) space} $\mathcal{H}_C$ is defined as the completion of $C(\mathbb{B}^*)$ under its induced norm, denoted by $\|\cdot\|_{\mathcal{H}_C}$, and is canonically identified with a continuously embedded Hilbert subspace of $\mathbb B$. The laws of $Z+h$ and $Z$ are mutually absolutely continuous if $h\in\mathcal H_C$, meaning that they assign probability zero to exactly the same measurable sets. If $h\notin\mathcal H_C$, the laws are mutually singular, meaning that there exists a measurable set to which one law assigns probability one and the other assigns probability zero. Thus, the CM space describes the admissible shifts of the Gaussian law, while $\mathbb B$ contains its sample paths \citep{van2008reproducing}.

\begin{example}\textup{\textbf{(Finite-Dimensional Hilbert Space)}}
\label{ex: finite CM}
Let $\phi_1,\ldots,\phi_m$ be continuous real-valued functions on $[0,1]$ satisfying
\[
\langle\phi_j,\phi_k\rangle_{L^2[0,1]}
=\int_0^1\phi_j(t)\phi_k(t)\,dt=\delta_{jk},
\qquad
\delta_{jk}=\begin{cases}1,&j=k,\\0,&j\ne k.\end{cases}
\]
Define $\mathcal H_m=\operatorname{span}\{\phi_1,\ldots,\phi_m\}$ with $L^2[0,1]$ inner product. Write $\phi(t)=(\phi_1(t),\ldots,\phi_m(t))^\top$. Every $f\in\mathcal H_m$ has a unique representation $f(t)=\theta^\top\phi(t)$, and for $g(t)=\eta^\top\phi(t)$, we have $\langle f,g\rangle_{\mathcal H_m}=\theta^\top\eta,\ \|f\|_{\mathcal H_m}=\|\theta\|_2$. Take $\mathbb B=\mathcal H_m$ and define the centered Gaussian random element as
\[
Z_m(t)=\sum_{k=1}^m\xi_k\phi_k(t),\qquad
\xi_1,\ldots,\xi_m\overset{\mathrm{iid}}\sim N(0,1).
\]
Its CM space is $\mathcal H_m$ with exactly the above inner product. In particular, adding $\sigma Z_m$ to any element in  $\mathcal H_m$ is equivalent to adding $N(0,\sigma^2I_m)$ noise to its corresponding orthonormal coefficient vector. The corresponding covariance operator $C_m$ for $Z_m$ is given in the proof of Corollary~\ref{col: finite-dimension functional GDP}.
\end{example}

\begin{definition}[Reproducing Kernel Hilbert Space]
\label{def: RKHS}
Let $\mathcal{T}$ be an index set and let $K:\mathcal{T}\times\mathcal{T}\rightarrow\mathbb{R}$ be a symmetric positive semidefinite kernel. The reproducing kernel Hilbert space (RKHS) associated with $K$, denoted by $\mathcal{H}_K$, is the completion of $\operatorname{span}\{K(t,\cdot):t\in\mathcal{T}\}$ with inner product
\[
\left\langle
\sum_j a_jK(s_j,\cdot),\sum_k b_kK(t_k,\cdot)
\right\rangle_{\mathcal H_K}
=\sum_{j,k}a_jb_kK(s_j,t_k),
\]
where zero-norm elements are identified. 
\end{definition}

The completion is naturally identified with a Hilbert space of functions on $\mathcal{T}$ through the reproducing property $f(t)=\langle f,K(t,\cdot)\rangle_{\mathcal H_K},t\in \mathcal{T}$. In particular, by the Cauchy--Schwarz inequality, $|f(t)| \leq \|f\|_{\mathcal{H}_K} \sqrt{K(t,t)}$. Consequently, if $\|f\|_{\mathcal{H}_K}\leq B$, then $|f(t)|\leq B\sqrt{K(t,t)}$. This bound provides a direct connection between the RKHS norm and the pointwise sensitivity of function-valued statistics.

\begin{example}\textup{\textbf{(RKHS and Gaussian Process)}}
\label{ex: RKHS CM}
Let $K$ be a continuous symmetric positive semidefinite kernel on $[0,1]^2$. A centered Gaussian process with covariance kernel $K$, denoted by $Z_K\sim\operatorname{GP}(0,K)$, is a collection of real-valued random variables $\{Z_K(t):t\in[0,1]\}$ such that, for every finite collection $t_1,\ldots,t_T$, $\bigl(Z_K(t_1),\ldots,Z_K(t_T)\bigr)^\top \sim N(0,K_T)$, where $K_T=[K(t_j,t_k)]_{j,k=1}^T$. In particular, $\mathbb E[Z_K(s)Z_K(t)]=K(s,t)$. Assume that $Z_K$ admits a continuous version. Then $Z_K$ can be viewed as a Gaussian random element in $\mathbb B=C[0,1]$, equipped with its supremum norm $\|f\|_\infty = \sup_{t\in [0,1]} |f(t)|$ and Borel $\sigma$-algebra. The CM space $\mathcal{H}_{C_K}$ corresponding to $Z_K$ coincides with the RKHS $\mathcal{H}_K$ associated with $K$; in particular, they have the same inner product and norm.
\end{example}

Since $K$ is continuous on the compact set $[0,1]^2$, the RKHS $\mathcal H_K$ embeds continuously into $C[0,1]$. The continuous-version assumption ensures that the Gaussian noise itself takes values in $C[0,1]$ almost surely, although its sample paths need not belong to $\mathcal H_K$. The continuity assumption holds for several commonly used kernels, for example, the Brownian motion kernel $K(s,t)=\min\{s,t\}$, and the Gaussian kernel
\begin{equation}
\label{eq: Gaussian kernel}
K(s,t) = \exp\Set{-\dfrac{(s-t)^2}{2\rho^2}},\ \rho>0,\ s,t\in [0,1].
\end{equation}
The covariance operator $C_K$ and the isometric CM-space identification are established in the proof of Corollary~\ref{col: RKHS GDP}. 

\subsection{Differential Privacy}
\label{subsec: Differential Privacy}
Let a dataset $X=(x_1,\ldots,x_n)$ denote a collection of $n$ records, each $x_i$ is an element in some set $\mathcal{X}$, and $\mathcal{X}^n$ denotes the collection of all possible data sets. For two data sets $X, X' \in \mathcal{X}^n$, we denote $X\sim X'$ if they differ in at most one entry, that is, their \emph{Hamming distance} $d_H(X,X') = \sum_{i=1}^n \mathbf{1}\{x_i\neq x_i'\}\leq 1$. We call such $X$ and $X'$ as \emph{adjacent data sets}. A \emph{mechanism} $\mathcal{M}: \mathcal{X}^n\to \mathcal{Y}$ is a randomized algorithm taking a data set $X$ as input and outputting values in a measurable space $\mathcal{Y}$. 

At a high level, differential privacy quantifies the difference between outputs of the mechanism when taking two adjacent data sets as inputs. We use the GDP framework throughout the paper, where the difference between distributions is measured through hypothesis testing. Let $P$ and $Q$ denote the probability distributions of $\mathcal{M}(X)$ and $\mathcal{M}(X')$, respectively. Consider a hypothesis test between: $H_0: Y\sim P$, $H_1: Y\sim Q$, where $Y$ is the observed output. Let $\psi:\mathcal{Y}\to [0,1]$ be a rejection rule. We use $\alpha_\psi = \mathbb{E}_P[\psi]$ to denote the Type I error, and $\beta_\psi = 1-\mathbb{E}_Q[\psi]$ to denote the Type II error, respectively. The difficulty of distinguishing between $P$ and $Q$ is characterized by the \emph{trade-off function} $T(P, Q):[0, 1] \to [0, 1]$ defined as
\[
T(P,Q)(\alpha)=\inf_{\psi}
\left\{
\beta_\psi:\alpha_\psi\le \alpha
\right\},
\]
where the infimum is taken over all measurable rejection rules $\psi$ \citep{dong2022gaussian}.

The trade-off function characterizes the minimal type II error as a function of the prescribed Type I error. A large trade-off function indicates that the two distributions $P,\,Q$ are difficult to distinguish and corresponds to a stronger privacy guarantee. For $\mu\ge 0$, define $G_\mu = T\bigl(N(0,1),N(\mu,1)\bigr)$ to be the \emph{Gaussian trade-off function} between unit-variance Gaussian with mean difference $\mu$. Explicitly, $G_\mu$ can be expressed as $G_\mu(\alpha) = \Phi(\Phi^{-1}(1-\alpha)-\mu)$, where $\Phi$ denotes the standard normal CDF. With the notion of Gaussian trade-off function, GDP is formally defined as follows.

\begin{definition}[Gaussian Differential Privacy (GDP), \citealp{dong2022gaussian}]
\label{def: mu-GDP}
A mechanism $\mathcal{M}$ satisfies $\mu$-GDP if for all adjacent data sets $X \sim X' \in \mathcal{X}^n$,
\[
T(\mathcal{M}(X),\mathcal{M}(X')) \ge G_\mu.
\]
\end{definition}

Intuitively, if a mechanism $\mathcal{M}$ satisfies $\mu$-GDP, then testing $H_0: Y\sim \mathcal{M}(X)$ versus $H_1: Y\sim \mathcal{M}(X')$ is at least as hard as testing $H_0: Y\sim N(0,1)$ versus $Y\sim N(\mu, 1)$. The parameter $\mu$ quantifies the distinguishability between the outputs of a mechanism on adjacent data sets: smaller values of $\mu$ correspond to stronger privacy guarantees, but with less data usability. 

\begin{definition}[Global Sensitivity in a Normed Space]
\label{def: Global sensitivity in a normed space}
Let $(\mathcal{B}, \|\cdot\|_{\mathcal{B}})$ be a normed vector space, and let $S: \mathcal{X}^n \to \mathcal{B}$ be a statistic computed from a dataset $X \in \mathcal{X}^n$. The global sensitivity of $S$ with respect to $\|\cdot\|_{\mathcal{B}}$ is defined as
\[
\Delta_{\mathcal{B}}(S)
=
\sup_{X\sim X'}
\|S(X)-S(X')\|_{\mathcal{B}},
\]
\end{definition}

The sensitivity measures the maximum change in the statistic in the norm of the output space when one individual's data is modified. In finite-dimensional settings when $\mathcal{B} = \R^d$, $\ell_1$- and $\ell_2$-sensitivities, denoted as $\Delta_1$ and $\Delta_2$, are two common choices. 

Definition~\ref{def: Global sensitivity in a normed space} extends naturally to function-valued statistics taking values in a function space, while it may require that the statistic resides in some specific space (e.g., CM space) to satisfy the desired GDP property. In particular, we say a measurable summary $S$ is \emph{compatible} with a centered Gaussian random element $Z\sim \mathcal{N}_{\mathbb{B}}(0,C)$, if $S(X)\in\mathcal{H}_C$ for any $X\in \mathcal{X}^n$. We denote $\Delta_{\mathcal{H}_C}(S) = \sup_{X\sim X'} \|S(X)-S(X')\|_{\mathcal{H}_C}$. Compatibility ensures that this norm is well-defined, whereas finiteness of the sensitivity is a separate requirement.

In practice, sensitivity may be unbounded when the data are unbounded, making it difficult to calibrate noise for privacy mechanisms. \emph{Clipping} is a standard technique in DP for limiting the influence of individual data values. The clipping of $x\in \R$ to the interval $[a,b]$ is defined as
\[
\operatorname{clip}_{[a,b]}(x) = \min\{b,\max\{a,x\}\}.
\]

By restricting values to a bounded interval, clipping can help control the sensitivity of statistical queries and facilitate the calibration of noise required to achieve certain DP guarantees.

\begin{proposition}[The Gaussian Mechanism \citep{dong2022gaussian}]
    Let $S: \mathcal{X}^n \to \R^d$ be a function with $\ell_2$ sensitivity, $\Delta_2(S) = \sup_{X\sim X'} \|S(X)-S(X')\|_2$.
    For any $S$, the Gaussian Mechanism releases
    \[
    \mathcal{M}(X)=S(X)+Z,
    \]
    where $Z\sim N(0,\sigma^2 I_d)$. Then $\mathcal{M}(X)$ satisfies $\mu$-GDP if $\sigma = \Delta_2(S)/\mu$.
\end{proposition}

In finite-dimensions, the Gaussian mechanism provides a fundamental approach to achieving GDP by adding Gaussian noise calibrated to the $\ell_2$-sensitivity of the underlying statistic. 

The GDP satisfies two important properties: \emph{composition} and \emph{post-processing} \citep{dong2022gaussian}. Suppose $k$ independent mechanisms $\mathcal{M}_1,\ldots,\mathcal{M}_k$ operate on the same dataset and satisfy $\mu_1$-GDP, $\ldots$, $\mu_k$-GDP, respectively. Their composition $\mathcal{M}=(\mathcal{M}_1,\ldots,\mathcal{M}_k)$ satisfies $\mu$-GDP, where the overall privacy parameter is given by $\mu=\sqrt{\sum_{i=1}^{k}\mu_i^2}$. This property characterizes how privacy loss accumulates across multiple private releases. Second, the post-processing property states that if a mechanism $\mathcal{M}:\mathcal{X}^n\to\mathcal{Y}$ satisfies $\mu$-GDP, then for any measurable function $g$ mapping $\mathcal{Y}$ to an arbitrary measurable output space $\mathcal{R}$, the transformed mechanism $\mathcal{M}'(X)=g(\mathcal{M}(X))$ also satisfies $\mu$-GDP. In other words, any data-independent transformation of the privatized output preserves the original privacy guarantee without incurring additional privacy loss. 

These properties are particularly useful for analyzing the PTF mechanisms, both when combining multiple private releases and when applying subsequent estimation procedures to privatized data.

\section{Problem Setup}
\label{sec: Problem Setup}
In this section, we first specify our functional data model with necessary assumptions, then we define the privacy and estimation operators and formalize the proposed two possible orderings of fitting and privatizing: \emph{Fit-Then-Protect} (FTP) and \emph{Protect-Then-Fit} (PTF) with algorithmic implementation in a general separable Hilbert space $\mathcal{H}$. Our goal is to release a private estimate of the sample mean function under a prescribed GDP guarantee. Finally we discuss the two metrics, IMSE and MSE, for evaluating the two procedures of FTP and PTF under the corresponding spaces for each metric.

\subsection{Functional Data Model}
\label{subsec: functional data model}
Throughout the paper, the observed data consist of $n$ independent functional trajectories $\mathbf{x}_i=(x_{i1},\ldots,x_{iT})^\top\in\mathcal{X}=\mathbb{R}^T$ at common points $0\leq t_1<t_2<\ldots<t_T\leq 1$. We denote $X=(\mathbf{x}_1,\ldots,\mathbf{x}_n)^\top\in\mathbb R^{n\times T}$ as the observed data matrix, and denote $\mathbf{t} = (t_1,\ldots,t_T)^\top$ as the evaluation grid. Let $f_1,\ldots,f_n\in\mathcal H$ denote the latent functional trajectories in a general Hilbert space $\mathcal H$, and let $\mathbf{f}_i = (f_i(t_1),\ldots,f_i(t_T))^\top$ denote the evaluation vector of $f_i$ at $\mathbf{t}$. Let $\epsilon_1,\ldots,\epsilon_n\in \mathbb {R}^T$ denote the random observation noise. The functional data are observed according to
\begin{equation}
\mathbf{x}_i = \mathbf{f}_i + \mathbf{\epsilon}_i,\quad
i=1,\ldots,n,
\label{eq:data_generation}
\end{equation}
with the following assumptions:
\begin{enumerate} 
\item[\textbf{A.1}] (Bounded functional norm) $\|f_i\|_{\mathcal H}\le B, i=1,\ldots,n$, for some unknown constant $B>0$.
\item[\textbf{A.2}] (Gaussian observation noise) Conditional on $f_1,\ldots,f_n$, $\epsilon_1,\ldots,\epsilon_n \overset{i.i.d.}{\sim} N(0,\sigma_0^2I_T)$ for a common $\sigma_0\geq 0$.
\end{enumerate}
Under Assumption \textbf{A.2}, when $\sigma_0$ is assumed to be non-zero, the procedures do not require knowledge of its true value. On the other hand, some results consider the case when $\sigma_0=0$, which denotes the degenerate distribution at zero, so $x_{ij}=f_i(t_j)$.

Additional assumptions specific to finite-dimensional Hilbert spaces and RKHS will be stated in the corresponding subsections. Our objective is to estimate the sample mean function $\bar{f}=\dfrac{1}{n}\sum_{i=1}^nf_i,$ based on $X$, while satisfying a prescribed $\mu$-GDP guarantee.

\subsection{Privacy–Estimation Operators and Orderings}
\label{subsec: Privacy–Estimation Operators and Orderings}
To estimate an individual function from its observed trajectory, we define a generic function estimation operator $\mathcal F:\mathcal X\rightarrow \mathcal H$, given by
\begin{equation}
\begin{aligned}
    \mathcal{F}(\mathbf{x})
    =
    \argmin_{f\in\mathcal H}
    \ 
    \mathcal{L}(f;\mathbf{x})
    \quad
    \T{subject to}
    \quad
    \|f\|_{\mathcal H}\le B',
\end{aligned}
\label{eq: F_operator}
\end{equation}
where $\mathcal L$ denotes a prescribed empirical loss and $B'$ is the fitting radius one can select. We assume that each fitting problem admits a minimizer and use a fixed deterministic measurable selection whenever the minimizer is not unique. In the RKHS setting, we select the minimizer with the smallest RKHS norm. For our functional data model, $\mathcal{L}$ also depends on $\mathbf{t}$, we omit this as the grid is treated as given throughout the paper. For example, weighted square loss
\begin{equation}
    \mathcal{L}(f;\mathbf{x})
    =
    (\mathbf{f}-\mathbf{x})^\top W (\mathbf{f}-\mathbf{x}),
    \quad
    \mathbf{f} = (f(t_1),\ldots,f(t_T))^\top
    \label{eq: weight least square}
\end{equation}
is often considered, which depends on the choice of the weight matrix $W$. A common choice is to set the weight matrix proportional to the inverse covariance matrix of $\mathbf{x}$, i.e., $W \propto \Sigma^{-1},\ \Sigma := \Cov(\mathbf{x})$. When $W = I$, the loss reduces to squared error.

Next, we denote the privacy mechanism as an operator. Let $\mathcal{P}_\mu^{\mathcal H}:\mathcal H\rightarrow \mathbb{B}$ denote a Gaussian perturbation operator. The fitting space \(\mathcal H\) and ambient output space \(\mathbb B\) play different roles and need not coincide. Its composition with a dataset summary satisfies $\mu$-GDP after calibration to the summary's CM sensitivity, and let $\mathcal{P}_\mu^{\mathcal{X}^n}:\mathcal{X}^n\rightarrow\mathcal{X}$ denote a $\mu$-GDP mechanism acting directly on the observed dataset. The latter may involve performing operations such as clipping, averaging, and Gaussian perturbation before releasing a privatized summary. The particular constructions of these mechanisms will be developed in Section~\ref{sec: Theoretical Results}.

Under the previous framework, the two privacy-estimation pipelines considered in this paper can be written compactly as follows.
\paragraph{Fit-Then-Protect (FTP).}
Each individual trajectory is first converted into a function through a functional estimator $\mathcal F_\mathrm{FTP}$, then estimated functions are averaged, and the resulting function is privatized:
\begin{equation}
\widetilde{f}_{\mathrm{FTP}}
=
\mathcal{P}_\mu^{\mathcal H}
\left(
\dfrac{1}{n}
\sum_{i=1}^n
\mathcal{F}_{\mathrm{FTP}}(\mathbf{x}_i)
\right).
\label{eq: FTP_operator}
\end{equation}

\paragraph{Protect-Then-Fit (PTF).}
Privacy is first enforced on the observed dataset, after which the function estimator $\mathcal F_\mathrm{PTF}$ is applied to the released summary:
\begin{equation}
\widetilde{f}_{\mathrm{PTF}}
=
\mathcal{F}_{\mathrm{PTF}}
\left(
\mathcal{P}_\mu^{\mathcal{X}^n}(X)
\right).
\label{eq: PTF_operator}
\end{equation}

Importantly, although both pipelines will satisfy the same prescribed privacy guarantee after appropriate calibration, they differ fundamentally in the order in which smoothing and privatization are performed. The primary goal of this paper is to understand how this ordering affects statistical accuracy.

\begin{remark}[Relation to Classical Smoothing Orderings]
In classical FDA, smoothing and aggregation \citep{ramsay2005functional, yao2005functional} play an analogous role as our ordering problem. One may either smooth each individual trajectory and then average, $\hat{f}_{\mathrm{STA}} = n^{-1} \sum_{i=1}^n \mathcal{F}(\mathbf{x}_i)$, or average the discrete observations and then smooth, $\hat{f}_{\mathrm{ATS}} = \mathcal{F}\bigl( n^{-1} \sum_{i=1}^n \mathbf{x}_i \bigr)$. We refer to these two procedures as \emph{Smooth-Then-Average} (STA) and \emph{Average-Then-Smooth} (ATS), respectively. When $\mathcal{F}$ is linear, the two estimators coincide, but they diverge when $\mathcal{F}$ involves active constraints or nonlinear regularization. Importantly, even in noiseless settings where STA and ATS produce identical estimates, the private pipelines FTP and PTF remain fundamentally different because privacy perturbations are injected at different stages of the estimation pipeline. We formally investigate this distinction in Section~\ref{subsec: Statistical Comparison}.
\end{remark}

In general, the two fitting operators $\mathcal{F}_{\mathrm{FTP}}$ and $\mathcal{F}_{\mathrm{PTF}}$ may differ as they have different inputs. The FTP estimator fits each individual trajectory before privatization, so its fitting loss can account for observation noise without needing to accommodate privacy noise. The PTF estimator operates on privacy-perturbed observations whose covariance structure may vary across sampling locations. Accordingly, $\mathcal{F}_{\mathrm{PTF}}$ may employ a loss function tailored to the privacy-induced noise, such as weighted least squares. When the privatized observations have the same covariance structure across sampling locations, the two operators may coincide.

\begin{algorithm}[t]
\caption{Fit-Then-Protect (FTP); $\mathtt{FTP}(X, \mathbf{t};B',\mu)$}
\label{alg: FTP}
\begin{algorithmic}[1]
\Statex \textbf{Fixed specification:} Gaussian law of $Z$ on $\mathbb B$, its CM space $\mathcal H_C$, and fitting loss $\mathcal L_{\mathrm{FTP}}$.
\State \textbf{Input:} Observed dataset $X=(\mathbf{x}_1,\ldots,\mathbf{x}_n)^\top$; public grid $\mathbf t=(t_1,\ldots,t_T)^\top$; constraint radius $B'>0$; privacy budget $\mu>0$.
\For{$i=1,\ldots,n$} 
\State \textbf{[Fitting]} $\hat{f}_i^{(c)} \gets 
\argmin_{f\in \mathcal{H}_C,\|f\|_{\mathcal{H}_C}\leq B'} 
\mathcal{L}_{\mathrm{FTP}}(f;\mathbf{x}_i)
$. 
\EndFor 
\State
$\overline{\hat {f}^{(c)}}\gets n^{-1}\sum_{i=1}^{n} \hat{f}_i^{(c)}$. 
\State 
Set sensitivity bound $\Delta\gets 2B'/n$ and $\sigma\gets\Delta/\mu$.\\
Sample $Z$ from its fixed Gaussian law, independently of $X$.\\
\textbf{[Protecting]}  
$\widetilde f_{\mathrm{FTP}} \gets
\overline{\hat {f}^{(c)}} + \sigma Z$. 
\State \textbf{Output:} $\widetilde f_{\mathrm{FTP}}$.
\end{algorithmic}
\end{algorithm}

\begin{algorithm}[t]
\caption{Protect-Then-Fit (PTF); $\mathtt{PTF}(X, \mathbf{t}; B',\mu, \mathbf{c})$}
\label{alg: PTF}
\begin{algorithmic}[1]
\Statex \textbf{Fixed specification:} Fitting space $\mathcal{H}$ and loss $\mathcal L_{\mathrm{PTF}}$.
\State \textbf{Input:} Observed dataset $X=(\mathbf{x}_1,\ldots,\mathbf{x}_n)^\top$; public grid $\mathbf t$; constraint radius $B'>0$; privacy budget $\mu>0$; clipping bounds $\mathbf c=(c_1,\ldots,c_T)^\top$, $c_j>0$.
\For{$j=1,\ldots, T$}
\State $\bar x_j^{(\mathrm{clip})}\gets n^{-1}\sum_{i=1}^n\operatorname{clip}_{[-c_j,c_j]}(x_{ij})$.
\State $\Delta_j\gets 2c_j/n$; $\sigma_{\mathrm{DP},j}\gets\Delta_j\sqrt T/\mu$.
\State Independently sample $Z_j\sim N(0,\sigma_{\mathrm{DP},j}^2)$, independently of $X$.
\State \textbf{[Protecting]} $\widetilde x_j\gets\bar x_j^{(\mathrm{clip})}+Z_j$.
\EndFor
\State $\widetilde{\mathbf x}\gets(\widetilde x_1,\ldots,\widetilde x_T)^\top$.
\State \textbf{[Fitting]} $\widetilde f_{\mathrm{PTF}}\gets\argmin_{f\in\mathcal H,\,\|f\|_{\mathcal H}\le B'}\mathcal L_{\mathrm{PTF}}(f;\widetilde{\mathbf x})$.
\State \textbf{Output:} $\widetilde f_{\mathrm{PTF}}$.
\end{algorithmic}
\end{algorithm}



While Equation~\eqref{eq: FTP_operator} and~\eqref{eq: PTF_operator} introduce FTP and PTF as operators, their implementation requires specifying the function space, fitting loss, fitting radius, clipping bounds, and privacy mechanism. Algorithms~\ref{alg: FTP} and~\ref{alg: PTF} provide concrete realizations of FTP and PTF, respectively. For FTP, we choose a $\mathbb{B}$-valued centered Gaussian random element $Z$ and fit individual trajectories in its Cameron--Martin space $\mathcal H_C$. Here, $\mathbb B\supset\mathcal{H}_C$ is the ambient space containing the Gaussian noise and the released estimator. In the finite-dimensional Hilbert space and RKHS settings considered later in Section~\ref{sec: Theoretical Results}, we use the Gaussian noise constructions in Examples~\ref{ex: finite CM} and~\ref{ex: RKHS CM}, respectively, so that $\mathcal H_C=\mathcal H$ with the same inner product and norm. The clipping constant $\mathbf{c} = (c_1,\ldots,c_T)^\top$ for PTF does not affect the DP guarantee, provided the privacy mechanism is calibrated accordingly, but can affect the estimation performance, which is similar to $B'$. In Section~\ref{subsec: Statistical Comparison}, we consider $\mathbf{c}$ that depends on both the function space under consideration and the sampling grid $\mathbf{t}$. For example, in the RKHS, the reproducing property implies $|f(t)|\le\|f\|_{\mathcal H_K}\sqrt{K(t,t)}$, so that under the Gaussian kernel, where $K(t,t)=1$, one may simply take $c_j=B'$. 

\subsection{Utility Criteria}
\label{subsec: Utility Criteria}
To compare the two procedures, we evaluate their estimation accuracy relative to the sample mean function $\bar f$. The choice of evaluation metric depends on the function space and the identifiability of the underlying function. In finite-dimensional function spaces, under suitable conditions, the underlying function can be recovered, allowing us to assess the global accuracy of the estimator relative to the target function. For finite-dimensional function spaces equipped with the $L^2[0,1]$ norm, we consider the integrated mean squared error (IMSE), where expectations are over observation and privacy noise conditional on $f_1,\ldots,f_n$ (similar for MSE).
\begin{equation}
\operatorname{IMSE}(\widetilde{f})
=
\int_0^1
\mathbb E
\left[
(\widetilde f(t)-\bar f(t))^2
\right]
dt.
\label{eq: IMSE}
\end{equation}

In infinite-dimensional function spaces, however, multiple functions may interpolate the same observed data points, making the underlying function potentially unidentifiable from the observed data alone. In this case, we evaluate estimation accuracy at the observed data points, and consider the empirical mean squared error (MSE) instead,
\begin{equation}
\operatorname{MSE}(\widetilde{f})
=
\frac1T
\sum_{j=1}^T
\mathbb E
\left[
(\widetilde f(t_j)-\bar f(t_j))^2
\right].
\label{eq: Empirical MSE}
\end{equation}

The remainder of the paper develops concrete realizations of the privacy operators in different function spaces and compares the statistical performance of the FTP and PTF pipelines under a common privacy budget. We write $\operatorname{IMSE}_{\mathrm{FTP}} = \operatorname{IMSE} (\widetilde{f}_{\mathrm{FTP}})$ and similar notation for brevity in the following sections.

\section{Theoretical Results}  \label{sec: Theoretical Results}
In this section, we develop the GDP mechanisms for functionals and compare the FTP and PTF regimes under the IMSE/MSE metric (depending on the function space). In Section~\ref{sec: GDP for functionals}, we extend the GDP mechanisms for functionals to the general separable Banach space and apply it to finite-dimensional Hilbert spaces and RKHSs. In Section~\ref{subsec: Statistical Comparison}, we give results for IMSE/MSE Ratio between the two regimes under the case without observation noise.
\subsection{Privacy Guarantees}
\label{sec: GDP for functionals}
Gaussian-process perturbations provide a natural way to privatize
function-valued summaries while preserving their functional structure.
\citet{hall2013differential} established $(\epsilon,\delta)$-DP guarantees
for releases evaluated at arbitrary finite collections of points, with
sensitivity measured in the RKHS associated with the noise covariance
kernel. \citet{mirshani2019formal} developed a Banach-space formulation
that provides privacy guarantees for full-function releases.
More recently, \citet{raisa2024noise} established a GDP guarantee for
Gaussian-process perturbations under RKHS.
The following theorem extends this GDP guarantee to Gaussian random
elements in a separable Banach space, with sensitivity measured in
the corresponding Cameron--Martin space.

\begin{theorem}[GDP for Function-Valued Summaries]
\label{thm: GDP_functional}
    Let $\mathbb{B}$ be a real separable Banach space, and let $Z\sim \mathcal{N}_{\mathbb{B}}(0,C)$ with the corresponding Cameron--Martin space $\mathcal{H}_C$. Let $S:\mathcal{X}^n \rightarrow \mathbb{B}$ be a measurable deterministic summary. Assume that $S$ is compatible with $Z$ and that $\Delta_{\mathcal{H}_C}(S)<\infty$. For any $\mu>0$ and $\sigma>0$ satisfying $\sigma\ge\Delta_{\mathcal{H}_C}(S)/\mu$, define
    \[
    \widetilde{S}(X)
    =
    S(X)+\sigma Z,
    \]
    Then, $\widetilde{S}(X)$ satisfies $\mu$-GDP.
\end{theorem}

\begin{proof}\textbf{sketch.}
Our proof works directly with the hypothesis-testing interpretation of GDP. Using the Cameron--Martin density formula, it reduces the comparison of adjacent output distributions to a real-valued Gaussian testing problem, yielding the privacy guarantee directly through its tradeoff function.    
\end{proof}

Theorem~\ref{thm: GDP_functional} provides the privacy guarantee for the FTP pipeline in Equation~\eqref{eq: FTP_operator}. In this setting, the summary is $S(X)=n^{-1}\sum_{i=1}^n\mathcal F(x_i)$. Since each fitted function belongs to the Cameron--Martin ball of radius $B'$, this summary is compatible with the Gaussian noise and has global sensitivity at most $2B'/n$. Consequently, adding $(2B'/(n\mu))Z$ ensures $\mu$-GDP. 

In general, the Cameron--Martin space $\mathcal H_C$ is a Hilbert subspace of $\mathbb B$ with potentially different norm. To apply Theorem~\ref{thm: GDP_functional} in the two settings considered below, we use the Gaussian noise constructions in Examples~\ref{ex: finite CM} and~\ref{ex: RKHS CM}. For the finite-dimensional space $\mathcal H_m$, we take $\mathbb B=\mathcal H_m$ and $Z=Z_m = \sum_{k=1}^m\xi_k\phi_k,\ \xi_1,\ldots,\xi_m\overset{\mathrm{iid}}\sim N(0,1)$. For the RKHS $\mathcal H_K$, we take $\mathbb B=C[0,1]$ and $Z=Z_K\sim\operatorname{GP}(0,K)$ with the continuous version specified in Example~\ref{ex: RKHS CM}. These constructions give $\mathcal H_{C}=\mathcal H_m$ and $\mathcal H_{C}=\mathcal H_K$, respectively, with equality of their inner products and norms. The following corollaries and the FTP--PTF comparisons in Section~\ref{subsec: Statistical Comparison} use these specified noise constructions, where the fitting space and corresponding CM space are equal.
\begin{corollary}\textup{\textbf{(Finite-Dimensional Hilbert Space Realization)}}
    Let $\mathcal H_m$ be a finite-dimensional Hilbert space, $\phi$ be the orthonormal basis, and $Z_m$ be the centered Gaussian random element as in Example~\ref{ex: finite CM}, and let $S(X)=\theta(X)^\top\phi$ be compatible with $Z_m$, where $\theta: \mathcal{X}^n\to \R^m$ is a measurable orthonormal coefficient summary. Suppose $\Delta_{\mathcal{H}_m}(S)<\infty$. For any $\mu>0$ and $\sigma>0$ with $\sigma\ge\Delta_{\mathcal{H}_m}(S)/\mu$, the release
    \[
    \widetilde S(X)=S(X)+\sigma Z_m=(\theta(X)+\sigma\xi)^\top\phi,
    \qquad \xi\sim N(0,I_m),
    \]
    is $\mu$-GDP. In particular, Algorithm~\ref{alg: FTP} with CM space $\mathcal H_m$ uses $\sigma=2B'/(n\mu)$ and adds Gaussian noise $\sigma\xi$ to the orthonormal coefficient vector.
    \label{col: finite-dimension functional GDP}
\end{corollary}

\begin{corollary}[RKHS Realization]
\label{col: RKHS GDP}
Let $K$ be a continuous symmetric positive semidefinite kernel on $[0,1]^2$ and a centered Gaussian process $Z_K\sim\operatorname{GP}(0,K)$ be as in Example~\ref{ex: RKHS CM}. For a measurable $S:\mathcal X^n\to\mathcal H_K$ that is compatible with $Z_K$, suppose $\Delta_{\mathcal{H}_K}(S)<\infty$. Then for any $\mu>0$, $\sigma>0$, and $\sigma\ge\Delta_{\mathcal{H}_K}(S)/\mu$, the release
\[
\widetilde S(X) = S(X)+\sigma Z_K
\] 
is $\mu$-GDP. Algorithm~\ref{alg: FTP} with CM space $\mathcal H_K$ therefore uses $\sigma=2B'/(n\mu)$. Its evaluations at the public grid satisfy $\widetilde{\mathbf f}_{\mathrm{FTP}} =\overline{\widehat{\mathbf f}^{(c)}}+(2B')/(n\mu)K_T^{1/2}\xi,\ \xi\sim N(0,I_T).$
\end{corollary}

Corollary~\ref{col: finite-dimension functional GDP} shows that the proposed functional GDP mechanism naturally extends the classical Gaussian mechanism from finite-dimensional vectors to finite-basis functional release through its coefficient vector. Corollary~\ref{col: RKHS GDP} specifies a full GP perturbation and its exact grid marginal. These constructions give direct implementations of Algorithm~\ref{alg: FTP}: the finite-dimensional Hilbert space release is represented by its coefficient vector, whereas an RKHS release evaluated only on the grid can be generated from the exact Gaussian grid marginal.

\begin{proposition}[Privacy Guarantee for PTF]
For fixed grid $\mathbf t$, fitting radius $B'>0$, privacy budget $\mu>0$, and clipping bounds $c_j>0$ for each $j$, the estimator $\widetilde{f}_{\mathrm{PTF}}$ given by 
Algorithm~\ref{alg: PTF} satisfies $\mu$-GDP.
\label{prop: GDP for PTF}
\end{proposition}
Proposition~\ref{prop: GDP for PTF} establishes the privacy guarantee for the PTF pipeline. The argument combines GDP composition for the privatized coordinate means with the post-processing property: once the privatized vector $\widetilde{x}$ satisfies $\mu$-GDP, applying the deterministic measurable fitting operator $\mathcal F_{\mathrm{PTF}}$ preserves this guarantee.

\subsection{Comparison of FTP and PTF}
\label{subsec: Statistical Comparison}
The privacy guarantees established in the previous subsection hold for a general Banach space. The accuracy of the FTP and PTF pipelines, however, depends strongly on the structure of the underlying estimation problem. In this subsection we compare the two procedures under two representative settings, finite-dimensional Hilbert space and RKHS. For each function space, we first analyze the error rates for FTP and PTF separately, then give the upper bound of the error ratio of FTP to PTF in the noiseless setting, which reveals the intrinsic effect of the ordering of estimation and privatization.

\subsubsection{Finite-dimensional Hilbert spaces}
We first study finite-dimensional Hilbert spaces, where we use the IMSE for the comparison. Consider the space $\mathcal H_m$ and continuous $L^2[0,1]$-orthonormal basis of Example~\ref{ex: finite CM}. Write $f_i(t)=\theta_i^\top\phi(t)$ and denote
\[
\Phi=[\phi_k(t_j)]_{j=1,\ldots,T;\,k=1,\ldots,m},\qquad M=\Phi^\top\Phi.
\]
with the following assumption
\begin{itemize}
    \item [\textbf{A.3}] $M = \Phi^\top\Phi$ is positive definite and $\|\phi(t_j)\|_2>0,\ j=1,\ldots,T$.
\end{itemize}
The requirement that $\|\phi(t_j)\|_2>0$ for each $j$ makes the pointwise clipping scales and inverse weights well defined. The positive definite condition in Assumption \textbf{A.3} depends on both the choice of orthonormal basis $\phi(t)$ and the evaluation points $t_1,\ldots,t_T$. In particular, $T\geq m$ is necessary. If $\phi(t)$ form an orthonormal basis of the polynomial space of degree at most $m-1$, then Assumption \textbf{A.3} holds if and only if the evaluation points $\mathbf{t}$ contain at least $m$ distinct values. For a general basis $\phi(t)$, the corresponding condition on the evaluation points may differ.

For a symmetric positive-definite matrix $A$, the \emph{spectral condition number} is defined as 
$\kappa(A)=\lambda_{\max}(A)/\lambda_{\min}(A)$, 
where $\lambda_{\min}(A)\ \text{and}\ \lambda_{\max}(A)$ denote the smallest and largest eigenvalues of $A$, respectively. The closer a condition number is to 1, the more uniform the eigenvalues are.

\begin{theorem}\textup{\textbf{(FTP IMSE Upper Bound in Finite-Dimensional Hilbert Spaces)}}
    Under the finite-dimensional Hilbert space $\mathcal H_m$ with Assumptions A.1, A.2, and A.3, consider the Gaussian perturbation $Z_m$ as in Example~\ref{ex: finite CM} with its CM space $\mathcal H_C=\mathcal H_m$. For the estimator $\widetilde{f}_{\mathrm{FTP}}$ given by $\mathtt{FTP}(X, \mathbf{t}; B',\mu)$ (Algorithm~\ref{alg: FTP}) with $\mathcal{L}$ defined in Equation~\eqref{eq: weight least square} and $W = I$. When $\sigma_0>0$ and $B' \geq B + \dfrac{\sqrt{m}\sigma_0}{\sqrt{\lambda_{\min}(M)}}$, $\widetilde{f}_{\mathrm{FTP}}$ satisfies 
    \begin{equation}
        \begin{aligned}
            \mathrm{IMSE}_{\mathrm{FTP}}
            \le& m\cdot\dfrac{4B'^2}{n^2\mu^2} + 
            \dfrac{\kappa(M)\sigma_0^2}{n}
            \operatorname{tr}
            (M^{-1})\\
            &+
            \Bra{1-\dfrac{1}{n}}
            \sqrt{3}\kappa(M)\sigma_0^2\operatorname{tr}
            (M^{-1})
            \exp
            \SqBra{
            -\dfrac{1}{4}\Bra{
            \dfrac{(B'-B)\sqrt{\lambda_{\min}(M)}}{\sigma_0}-\sqrt{m}
            }^2
            }.
        \end{aligned}
    \end{equation}
    When $\sigma_0=0$, for $B'\ge B$ the exact error is $\operatorname{IMSE}_{\mathrm{FTP}}=4mB'^2/(n^2\mu^2)$.
    \label{thm: FTP IMSE}
\end{theorem}
Theorem~\ref{thm: FTP IMSE} gives the upper bound for IMSE of FTP under finite-dimensional Hilbert spaces. The three terms correspond to privacy noise, observation error, and the effect of an active constraint. Fix the basis $\phi$, grid $\mathbf{t}$, $B$, and $\sigma_0>0$. Strengthening the lower bound on $B'$ in
Theorem~\ref{thm: FTP IMSE} by replacing $\sqrt m$ with $\sqrt m+a\sqrt{\log n}$, for some $a>2$, the constraint remainder becomes $O(n^{-a^2/4})=o(n^{-1})$, yielding
\[
\operatorname{IMSE}_{\mathrm{FTP}}
=
O\!\left(\frac1n+\frac{(B')^2}{n^2\mu^2}\right).
\]
For a fixed $\mu>0$, if we additionally have $B'=O(\sqrt n)$, it suffices to obtain $O(n^{-1})$ IMSE. 

The dependence on the sampling grid can also be understood through the observation-error term. For a fixed continuous orthonormal basis, if the empirical measures of a sequence of grids converge weakly to the uniform distribution on $[0,1]$ as $T\to\infty$, then $M/T\to I_m$, $\kappa(M)\to1$, and $\operatorname{tr}(M^{-1})\sim m/T$. Consequently, the observation-error term in the theorem is asymptotically $\sigma_0^2m/(nT)$.

\begin{theorem}\textup{\textbf{(PTF IMSE Lower Bound in Finite-Dimensional Hilbert Spaces)}}
\label{thm: PTF IMSE lower}
Under the finite-dimensional Hilbert space $\mathcal{H}_m$ with Assumptions A.1, A.2, and A.3. Denote $\tau=\dfrac{2B'\sqrt{T}}{n\mu}.$ For the estimator $\widetilde{f}_{\mathrm{PTF}}$ given by $\mathtt{PTF}(X, \mathbf{t}; B',\mu, \mathbf{c})$ (Algorithm~\ref{alg: PTF}) with $\mathcal{L}$ defined in Equation~\eqref{eq: weight least square}, $W=\operatorname{diag}\{\|\phi(t_j)\|_2^{-2}:j=1,\ldots,T\}$, $c_j = B'\|\phi(t_j)\|_2$. Denote $A = \Phi^\top W\Phi$, then $A$ is positive definite. When $\sigma_0 = 0$ and $B'\geq B + \dfrac{\tau\sqrt m}{\sqrt{\lambda_{\min}(A)}}$, $\widetilde{f}_{\mathrm{PTF}}$ satisfies
\begin{equation}
\begin{aligned}
\operatorname{IMSE}_{\mathrm{PTF}}
\geq
\tau^2\operatorname{tr}(A^{-1})
\Set{1-
\sqrt{3}
\exp\left[
-\frac14
\left\{
\frac{(B'-B)\sqrt{\lambda_{\min}(A)}}{\tau}
-\sqrt m
\right\}^2
\right]}.
\end{aligned}
\label{eq: PTF lower bound}
\end{equation}
\end{theorem}
Theorem~\ref{thm: PTF IMSE lower} gives the lower bound for IMSE of PTF under the finite-dimensional Hilbert space without observation noise, when $B'>B$ for sufficiently large $n$, the exponential term decays rapidly, which leads to our following corollary.
\begin{corollary}\textup{\textbf{(IMSE Ratio Upper Bound for Finite-Dimensional Hilbert Spaces)}}
    Consider the FTP and PTF estimators defined in Theorem~\ref{thm: FTP IMSE} and Theorem~\ref{thm: PTF IMSE lower} with $\sigma_0 = 0$. For any fixed $T, m, \mu$, and a non-decreasing sequence $B'_n$ such that $B'_n>B$ for sufficiently large $n$, we have
    \begin{equation}
        \begin{aligned}
            \limsup_{n\to\infty}
            \frac{\operatorname{IMSE}(\widetilde f_{\mathrm{FTP}})}
            {\operatorname{IMSE}(\widetilde f_{\mathrm{PTF}})}
            \le\frac{m}{T\operatorname{tr}(A^{-1})}\le\frac1m,
        \end{aligned}
    \end{equation}
    where $A=\Phi^\top W\Phi$ is defined in Theorem~\ref{thm: PTF IMSE lower}, and $B'_n$ serves the role of $B'$ in Theorems~\ref{thm: FTP IMSE} and~\ref{thm: PTF IMSE lower}.
    \label{col: IMSE Ratio for Polynomials}
\end{corollary}
Corollary~\ref{col: IMSE Ratio for Polynomials} suggests that for $m>1$, FTP has strictly smaller IMSE for all sufficiently large $n$ in this noiseless regime in finite-dimensional Hilbert spaces. The asymptotic upper bound $1/m$ on the FTP-to-PTF IMSE ratio decreases with the number of basis functions $m$. In particular, for the polynomial space, $m$ is distinct from the polynomial degree. When $m=2$, one may consider the orthonormal basis $\{1,\sqrt{3}(2t-1)\}$ for $\operatorname{span}\{1,t\}$, or the orthonormal basis $\{1,\sqrt{5}(3t^2-1)/2\}$ for $\operatorname{span}\{1,t^2\}$. Both choices give the asymptotic ratio bound $1/2$.
\subsubsection{RKHS}
Let $\mathcal{H}_K$ be the RKHS associated with kernel $K(\cdot,\cdot)$, denote
$K_T = \SqBra{K(t_k,t_l)}_{k,l = 1}^T$ 
with the following assumption:
\begin{itemize}
    \item [\textbf{A.4}] $K_T$ is positive definite and $K(t_j,t_j)>0,\ j=1,\ldots,T$.
\end{itemize}
The comparison is performed under the empirical MSE criterion
\begin{theorem}\textup{\textbf{(MSE Upper Bound for FTP in RKHS)}}
    Under the RKHS $\mathcal{H}_K$ with Assumptions A.1, A.2, and A.4, consider the Gaussian perturbation $Z_K\sim \operatorname{GP}(0,K)$ with its CM space $\mathcal{H}_C = \mathcal{H}_K$ be as in Example~\ref{ex: RKHS CM}. For the estimator $\widetilde{f}_{\mathrm{FTP}}$ given by $\mathtt{FTP}(X, \mathbf{t}; B',\mu)$ (Algorithm~\ref{alg: FTP}) with $\mathcal{L}$ defined in Equation~\eqref{eq: weight least square} and $W = I$. When $\sigma_0>0$ and $B'\geq B + \dfrac{\sigma_0\sqrt{T}}{\sqrt{\lambda_{\min}(K_T)}}$, we have
    \begin{align*}
        \mathrm{MSE}_{\mathrm{FTP}}
        &\leq \dfrac{4{B'}^2}{n^2\mu^2 T} \operatorname{tr}(K_T)
        + \dfrac{\sigma_0^2}{n}\\
        &\qquad + \left(1-\frac1n\right)
        \sigma_0^2
        \sqrt{\frac{T+2}{T}}
        \exp\left\{
        -\frac14
        \left(
        \frac{(B'-B)\sqrt{\lambda_{\min}(K_T)}}{\sigma_0}
        -\sqrt{T}
        \right)^2
        \right\}.
    \end{align*}
    When $\sigma_0=0$, for $B'\ge B$ the exact error is $\operatorname{MSE}_{\mathrm{FTP}}=4{B'}^2\operatorname{tr}(K_T)/(n^2\mu^2T)$.
    \label{thm: FTP MSE RKHS}
\end{theorem}
Theorem~\ref{thm: FTP MSE RKHS} provides a general upper bound for the MSE of FTP in the RKHS setting. The result is analogous to Theorem~\ref{thm: FTP IMSE}, but the present formulation allows the true functions $f_i$ to belong to the full, potentially infinite-dimensional RKHS. Fix the kernel $K$, grid $\mathbf t$, $B$, and $\sigma_0>0$, if the lower bound on $B'$ in
Theorem~\ref{thm: FTP MSE RKHS} is strengthened by replacing $\sqrt T$
with $\sqrt T+a\sqrt{\log n}$, for some $a>2$, then Theorem~\ref{thm: FTP MSE RKHS} gives
\[
\operatorname{MSE}_{\mathrm{FTP}}
=
O\!\left(\frac1n+\frac{(B')^2}{n^2\mu^2}\right).\]
For fixed $\mu>0$, the strengthened lower bound on $B'$, together with $B'=O(\sqrt n)$, ensures an $O(1/n)$ MSE. These conditions balance the growth of $B'$ needed to
control the constraint effect against the resulting increase in privacy noise. For the Gaussian kernel, $K(t_j,t_j)=1$ for every $j$, so $\operatorname{tr}(K_T)=T$ and the privacy contribution is exactly $4(B')^2/(n^2\mu^2)$.

\begin{theorem}[MSE Lower Bound for PTF in RKHS]
    \label{thm: PTF MSE lower bound}
    Under the RKHS $\mathcal{H}_K$ with Assumptions A.1, A.2, and A.4. Denote $\tau=\dfrac{2B'\sqrt{T}}{n\mu}.$ For the estimator $\widetilde{f}_{\mathrm{PTF}}$ given by $\mathtt{PTF}(X, \mathbf{t}; B',\mu, \mathbf{c})$ (Algorithm~\ref{alg: PTF}) with $\mathcal{L}$ defined in Equation~\eqref{eq: weight least square}, $W=\operatorname{diag}\{K(t_j,t_j)^{-1}:j=1,\ldots,T\}$, $c_j = B'\sqrt{K(t_j,t_j)}$. Denote $A = K_T^{1/2} W K_T^{1/2}$, then $A$ is positive definite.
    When $\sigma_0=0$ and $B'\geq B + \dfrac{\tau\sqrt T}{\sqrt{\lambda_{\min}(A)}}$, we have
    \begin{equation}
    \begin{aligned}
    \mathrm{MSE}_{\mathrm{PTF}}
    \geq
    \dfrac{\tau^2}{T}\operatorname{tr}(W^{-1})
    \SqBra{1
    -\sqrt{3}
    \exp\Set{
    -\frac14
    \Bra{
    \frac{(B'-B)\sqrt{\lambda_{\min}(A)}}{\tau}
    -\sqrt T
    }^2
    }}.
    \end{aligned}
    \label{eq: PTF MSE lower bound}
    \end{equation}
\end{theorem}
Theorem~\ref{thm: PTF MSE lower bound} provides a lower bound for MSE of PTF under RKHS without observation noise, which is similar to Theorem~\ref{thm: PTF IMSE lower}. In particular, under the Gaussian Kernel, we have $K(t_j,t_j)=1$ for each $j$, therefore $\operatorname{tr}(W^{-1})=T$.
\begin{corollary}[MSE Comparison under RKHS]
    Consider the FTP and PTF estimators defined in Theorem~\ref{thm: FTP MSE RKHS} and Theorem~\ref{thm: PTF MSE lower bound} with $\sigma_0 = 0$. Then, for any fixed $T, \mu$, and a non-decreasing sequence $B'_n$ such that $B'_n>B$ for sufficiently large $n$, we have
    \begin{equation}
        \begin{aligned}
            \limsup_{n\to\infty}
            \frac{\operatorname{MSE}(\widetilde f_{\mathrm{FTP}})}
            {\operatorname{MSE}(\widetilde f_{\mathrm{PTF}})}
            \leq \dfrac{1}{T}.
        \end{aligned}
    \end{equation}
    where $B'_n$ serves the role of $B'$ in Theorems~\ref{thm: FTP MSE RKHS} and~\ref{thm: PTF MSE lower bound}.
    \label{col: MSE ratio in RKHS}
\end{corollary}
Corollary~\ref{col: MSE ratio in RKHS} shows that in the noiseless setting with $T>1$, FTP achieves strictly smaller empirical MSE than PTF for sufficiently large $n$, with an asymptotic FTP-to-PTF MSE ratio bounded above by $1/T$. The comparison uses $\operatorname{tr}(W^{-1})=\operatorname{tr}(K_T)$, which applies to the kernels satisfying the continuous conditions in Example~\ref{ex: RKHS CM}.
\section{Simulation Studies}
\label{sec: Simulation Study}
In this section, we compare the empirical performance of the proposed strategies, \emph{Protect-Then-Fit} versus \emph{Fit-Then-Protect}, implemented through
Algorithms~\ref{alg: FTP} and~\ref{alg: PTF}. We focus on the RKHS associated with the Gaussian kernel; the corresponding simulation setup and results for polynomial spaces are presented in Appendix~\ref{sec: Additional Simulation Results}. In Section~\ref{subsec: Effect of Fitting Dimension}, we first examine the noiseless error ratios in relation to the theoretical bounds in Section~\ref{subsec: Statistical Comparison}, and in Section~\ref{subsec: Choices of Constraint Radius}, we investigate how different choices of the constraint radius $B'$ affect the performance of FTP and PTF. Finally, in Section~\ref{subsec: Effect of Observation Noise}, we further explore the effect of different magnitudes of observation noise. The RKHS simulation setup is specified below.

\textbf{RKHS Simulation Setup.} Consider the RKHS $\mathcal{H}_K$ induced by the Gaussian kernel $K$ defined in Equation~\eqref{eq: Gaussian kernel} with bandwidth $\rho=0.2$. The functions $f_1,\ldots,f_n\in\mathcal{H}_K$ are generated from $M=20$ kernel centers $\Set{c_1,\cdots,c_M}$ such that $f_i(t) = \alpha_i^\top K_c(t), K_c(t) = (K(c_1,t),\ldots,K(c_M,t))^\top$. We sample the coefficient vector $\alpha_i$ i.i.d. with a nonzero mean direction from a common bounded distribution, subject to $\alpha_i^\top K_{C}\alpha_i\le B^2$, where $K_{C}=[K(c_j,c_k)]_{j,k=1}^M$ and $B=1$. This ensures that $\|f_i\|_{\mathcal H_K}\le B$ while producing a nontrivial population mean function. The data matrix $X$ is observed based on the model: $x_{ij} = f_i(t_j)+\epsilon_{ij}, \epsilon_{ij}\overset{i.i.d.}{\sim} N(0,\sigma_0^2)$, with an equally spaced grid $t_j=(j-1)/(T-1), j=1,\ldots,T$. The latent functions are generated once for each latent-data setting and held fixed across repetitions. Each repetition uses fresh observation noise and privacy noise, with FTP and PTF applied to the same observed dataset. We estimate the conditional MSEs by averaging squared errors over Monte Carlo repetitions, holding the latent functions fixed. The reported error ratios are ratios of these Monte Carlo mean squared errors.

The constrained RKHS fitting in Algorithms~\ref{alg: FTP} and~\ref{alg: PTF} is implemented through finite kernel representations at the observation grid. An eigendecomposition of $K_T$ reduces each fitting problem to a scalar equation for the Lagrange multiplier associated with the RKHS norm constraint, which is solved numerically when the constraint is active. The implementation imposes a regularization floor of $10^{-8}$ for numerical stability. The FTP privacy perturbation is sampled on the observation grid from $N\!\left(0,(2B'/(n\mu))^2K_T\right)$.

\begin{figure}[t]
    \centering
    \includegraphics[width=\linewidth]{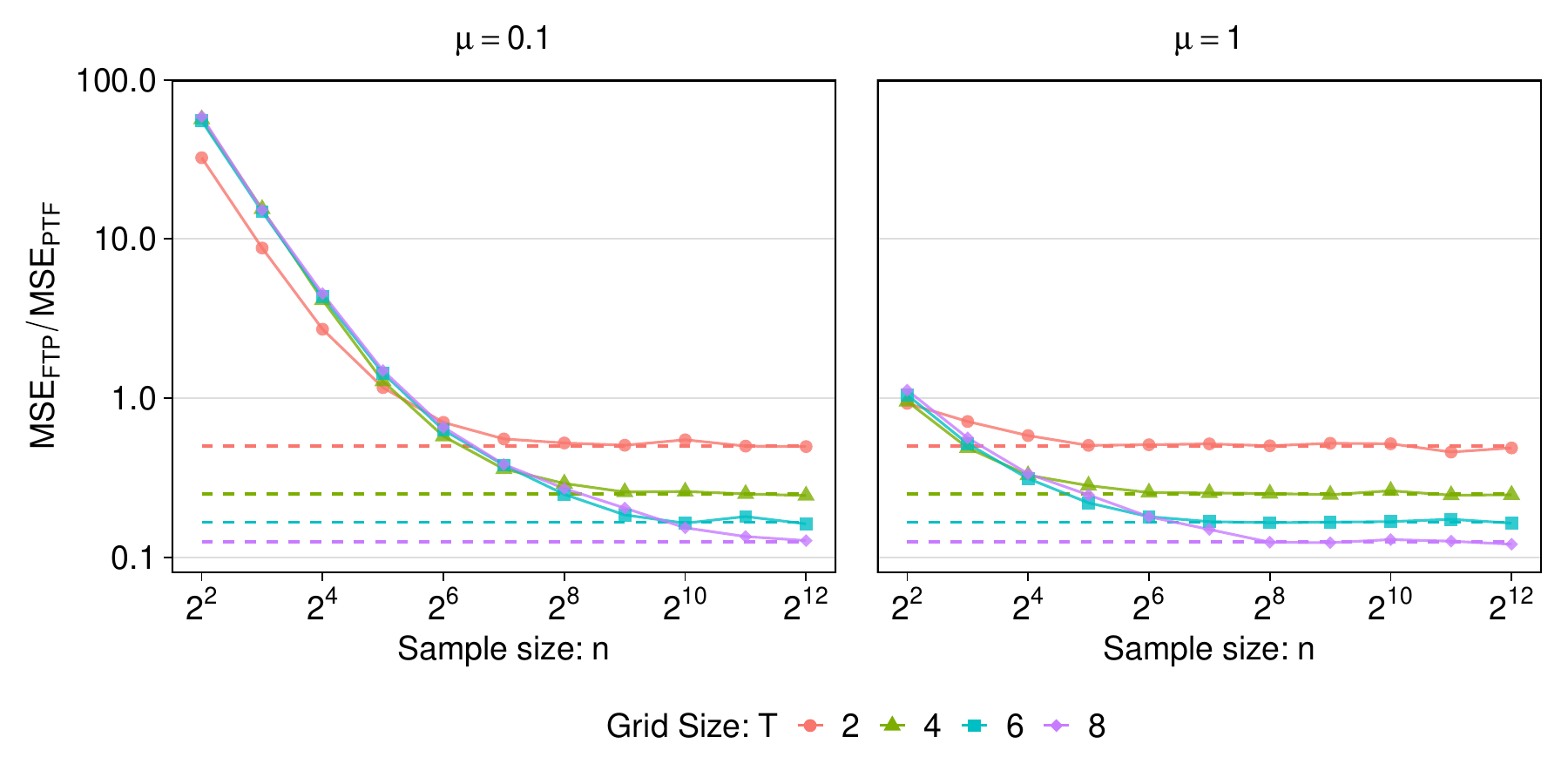}
    \caption{Empirical MSE ratio of FTP to PTF versus sample size $n$ in the noiseless RKHS setting, with $B'=1.1$ and $\mu\in\{0.1,1\}$. Colors and symbols identify the grid size $T$. Dashed lines show the corresponding $1/T$ asymptotic upper bound. Both axes use logarithmic scales.}
    \label{fig: MSE Ratio with Bound in RKHS}
\end{figure}

We vary the sample sizes, number of grid points, privacy budgets, constraint radius, and observation-noise standard deviations as follows: $n\in \Set{2^2, 2^3, \ldots, 2^{12}}$, $T\in \Set{2,4,6,8}$, $\mu\in \Set{0.1, 1}$, $B'\in\Set{B, 1.1 B, \sqrt{\log n}, \log n, \sqrt{n}, n^{2/3}}$, and $\sigma_0\in\Set{0, 0.1, 1, 10}$. The combinations used in each experiment are specified in the corresponding subsection. 

All the simulations are run with 1000 Monte Carlo repetitions. Simulation settings are evaluated in parallel in R using \texttt{foreach} and \texttt{doParallel}, with separate random-number streams assigned to the workers. 

\subsection{Effect of Fitting Dimension}
\label{subsec: Effect of Fitting Dimension}
We first investigate how the fitting dimension affects the relative accuracy of FTP and PTF. Here, the fitting dimension refers to the dimension of the space used to represent the fitted functions $\hat{f}_i^{(c)}, i=1,\ldots,n$ and $\widetilde{f}_{\mathrm{PTF}}$ in Algorithm~\ref{alg: FTP} and~\ref{alg: PTF}. In the RKHS setting, we use the minimum-norm constrained fits in $\operatorname{span}\{K(t_1,\cdot),\ldots,K(t_T,\cdot)\}$, which has dimension $T$ when $K_T$ is positive definite. Thus, varying the grid size $T$ also varies the fitting dimension, although the underlying RKHS is infinite-dimensional. We consider $T\in\{2,4,6,8\}$ and $n\in\{2^2,\ldots,2^{12}\}$, and set $\sigma_0=0$ to examine the noiseless comparison established in Corollary~\ref{col: MSE ratio in RKHS}. The privacy budgets $\mu\in\{0.1,1\}$ represent stronger and weaker privacy protection, respectively. We choose $B'=1.1 > B = 1$ as a representative case satisfying the corollary's condition. The dashed lines at $1/T$ provide baselines motivated by Corollary~\ref{col: MSE ratio in RKHS}.

Figure~\ref{fig: MSE Ratio with Bound in RKHS} shows that the MSE ratios of FTP to PTF generally decrease as sample size $n$ increases and approach the corresponding $1/T$ reference levels, with small Monte Carlo fluctuations. These findings are consistent with Corollary~\ref{col: MSE ratio in RKHS}, and suggest that FTP will outperform PTF asymptotically when $T\ge 2$. For the stronger privacy setting $\mu=0.1$, the ratios are substantially above one at small sample sizes, showing that PTF can be preferable in this regime. For larger $T$ or smaller $\mu$, a larger sample size is generally needed for the MSE ratio to approach the noiseless reference level $1/T$. Both changes increase the variance of the privacy noise in PTF and tend to make its fitting constraint active more frequently at a given $n$, as illustrated in Table~\ref{tab: RKHS noiseless details} in Appendix~\ref{sec: Additional Simulation Results}. The corresponding polynomial-space results are presented in Figure~\ref{fig: IMSE Ratio with Bound in Polynomial} in the same appendix.

\subsection{Choices of Constraint Radius}
\label{subsec: Choices of Constraint Radius}
We next examine the role of the fitting radius $B'$, which controls both the fitting constraint and the magnitude of the privacy noise. A smaller radius imposes stronger restrictions on the fitted functions, whereas a larger radius increases the noise required for privacy. We fix $\sigma_0=1$, $\mu=1$, and $T=4$, and compare $B'\in\{B,\sqrt{\log n},\log n,\sqrt n,n^{2/3}\}$ over $n\in\{2^2,\ldots,2^{12}\}$. The choice $B'=B$ serves as an oracle baseline as it represents the knowledge of the true functional norm bound. We include $\sqrt{\log n}$, $\log n$ and $\sqrt{n}$ as three representative rates based on our theoretical analysis in Section~\ref{subsec: Statistical Comparison}, and $B' = n^{2/3}$ is out of the conditions for an $O(1/n)$ MSE rate of FTP, whose FTP privacy contribution is of order $n^{-2/3}$. We plot the scaled MSE (multiplied by $n$) to help capture its rate, for example, an error of order $n^{-1}$ corresponds to a scaled error of constant order. 

\begin{figure}[t]
    \centering
    \includegraphics[width=\linewidth]{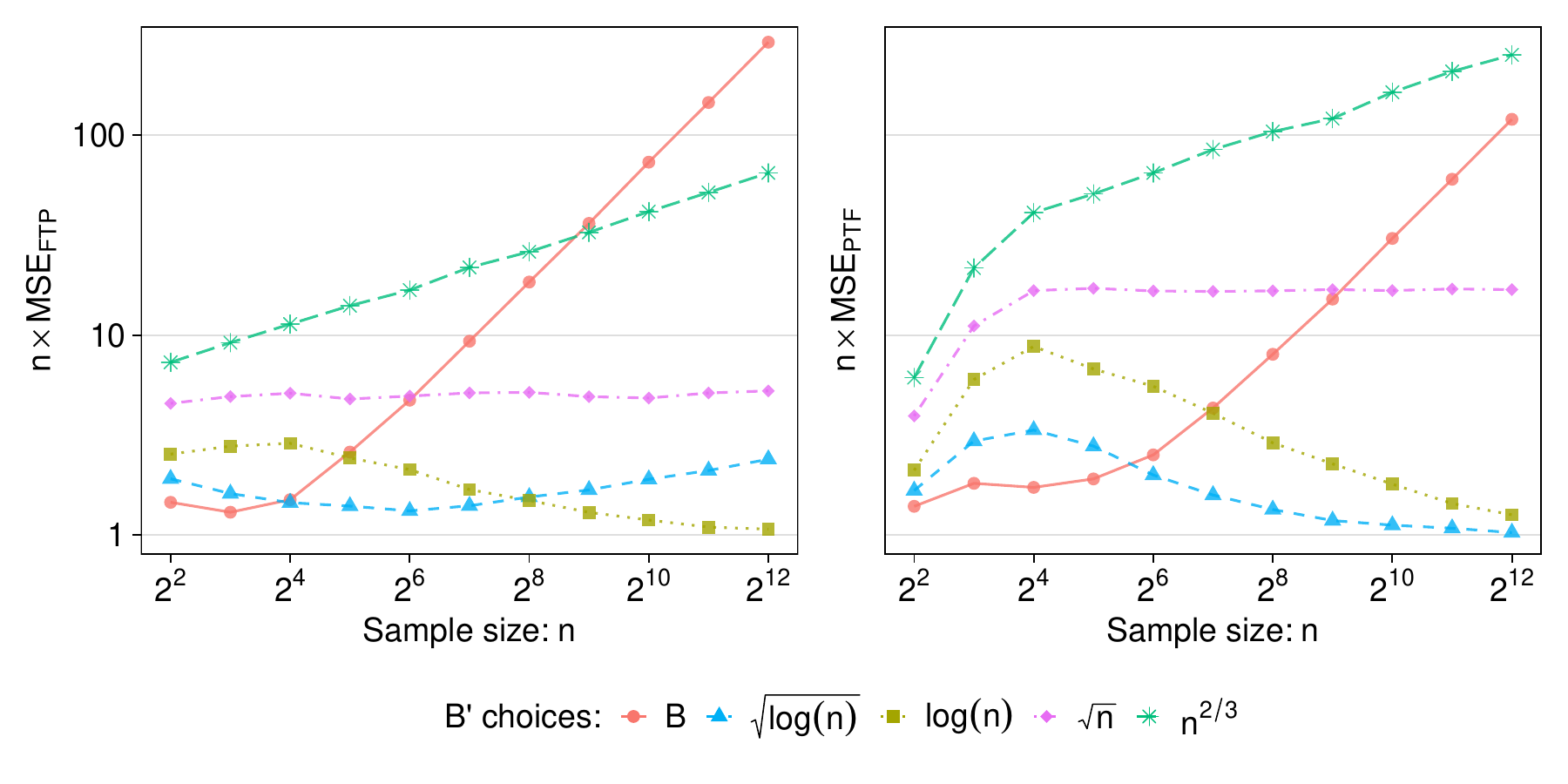}
    \caption{Scaled empirical MSE ($n\times\operatorname{MSE}$) for FTP (left) and PTF (right) under different fitting radii $B'$, with $\sigma_0=1$, $\mu=1$, and $T=4$. Both axes use logarithmic scales.}
    \label{fig: RKHS_Bprime_choices}
\end{figure}

Figure~\ref{fig: RKHS_Bprime_choices} shows that no single radius choice performs best across all sample sizes. The oracle radius $B'=B$ performs well at small sample sizes, but its scaled MSE increases substantially as $n$ grows in both pipelines. For FTP, $B'=\log n$ gives the smallest error among the displayed choices at the largest sample sizes, with $n\times\operatorname{MSE}_{\mathrm{FTP}}$ approaching 1. The choice $B'=\sqrt{\log n}$ initially performs well but exhibits a gradual increase in scaled error later as $n$ grows. For PTF, both logarithmic choices yield improved performance at larger sample sizes, with $\sqrt{\log n}$ achieving the smaller error near the end of the range. For both FTP and PTF, the $\sqrt n$ curves are approximately flat at larger sample sizes, whereas the $n^{2/3}$ curves exhibit an increasing trend. These patterns illustrate the trade-off between imposing overly restrictive constraints on the fitted functions and allowing the privacy-noise contribution to decay too slowly. In practice, one should select an appropriate fitting radius to balance the constraint effect against privacy noise. The polynomial-space results are reported in Figure~\ref{fig: Poly_Bprime_choices} in Appendix~\ref{sec: Additional Simulation Results}.

\subsection{Effect of Observation Noise}
\label{subsec: Effect of Observation Noise}
Finally, we investigate how observation noise affects the relative performance of FTP and PTF by comparing their MSE ratio under noisy observations. We vary $n\in\{2^2,\ldots,2^{12}\}$ and $T\in\{2,4,6,8\}$, and consider $\sigma_0\in\{0.1,1,10\}$. Based on the results and discussion in Section~\ref{subsec: Effect of Fitting Dimension} and~\ref{subsec: Choices of Constraint Radius}, we fix $\mu=1$ and $B'=\log n$. Since $B=1$ and the Gaussian kernel has unit diagonal, the latent function values satisfy $|f_i(t)|\le1$. Thus, the selected noise standard deviations range from relatively small to substantially larger than this bound. The dashed $1/T$ lines are retained as noiseless theoretical baselines.

\begin{figure}[t]
    \centering
    \includegraphics[width=\linewidth]{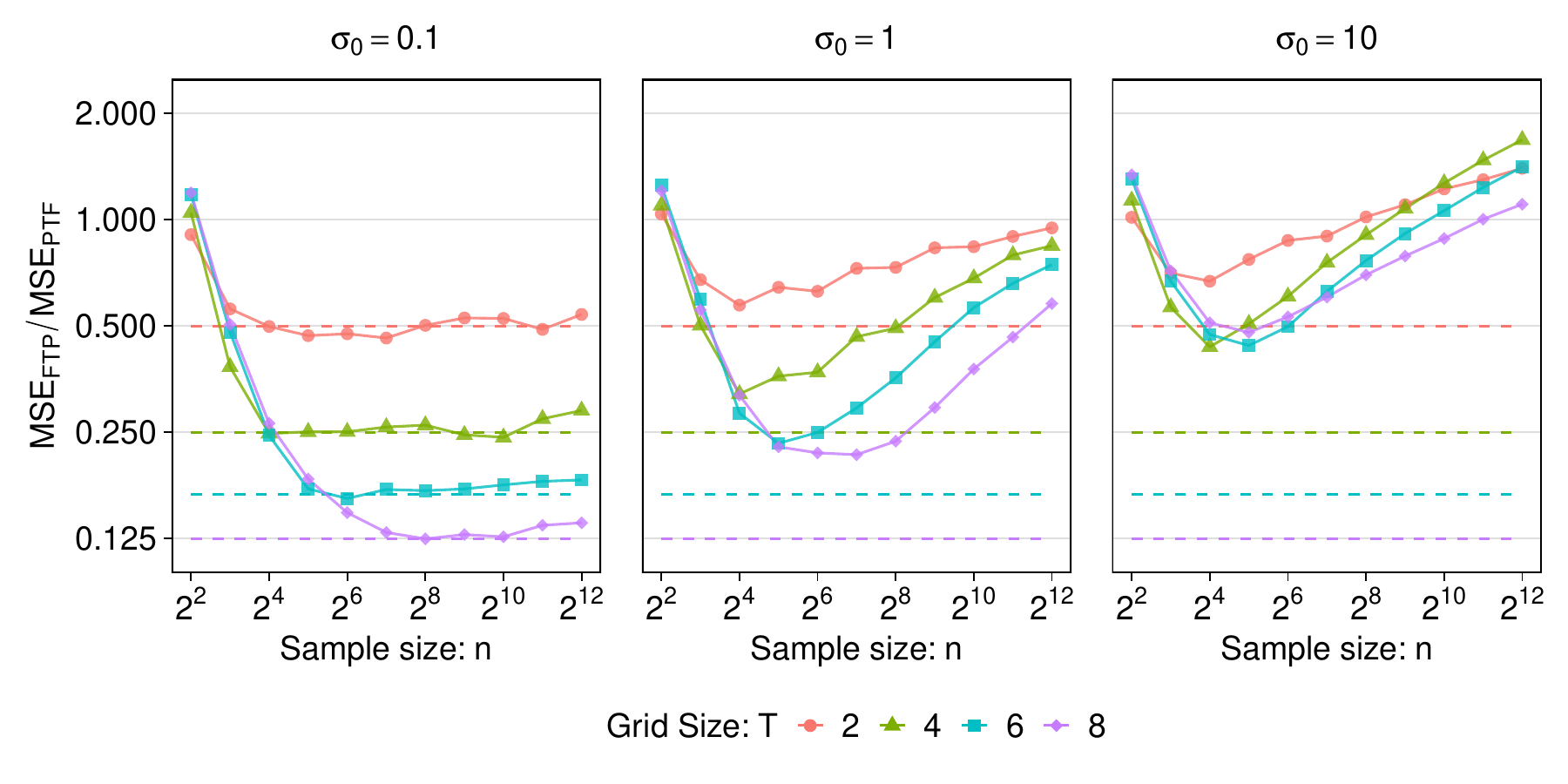}
    \caption{Empirical MSE ratio of FTP to PTF versus sample size $n$ for $\sigma_0\in\{0.1,1,10\}$, with $\mu=1$ and $B'=\log n$. Colors and symbols identify the grid size $T$. Dashed lines show the noiseless asymptotic reference bounds $1/T$. Both axes use logarithmic scales.}
    \label{fig: MSE Ratio with Bound in RKHS with noise}
\end{figure}

Figure~\ref{fig: MSE Ratio with Bound in RKHS with noise} shows that the relative performance of FTP and PTF depends strongly on the observation-noise level. When $\sigma_0=0.1$, the ratios decrease rapidly and remain close to the noiseless reference baselines, with FTP retaining a substantial advantage at larger sample sizes. When $\sigma_0=1$, the ratios initially decrease but then increase, indicating that FTP's relative advantage weakens at larger sample sizes, although all displayed ratios remain below one at the largest $n$. When $\sigma_0=10$, the ratios exhibit a brief initial decrease followed by a clear increasing trend, with all ratios exceeding one at the largest sample size, indicating that PTF achieves lower MSE in this regime. Overall, these findings demonstrate that FTP's advantage persists under low observation noise at sufficiently large sample sizes, extending the empirical evidence beyond the noiseless setting considered in our theory. The high-noise results further reveal the limitations of this advantage, showing that PTF can become preferable when observation noise is substantial. The corresponding results for polynomial-space are presented in Figure~\ref{fig: IMSE Ratio with Bound in Polynomial with noise} in Appendix~\ref{sec: Additional Simulation Results}.

\section{Discussion and Directions for Future Work}
\label{sec: Discussion}
In this paper, we study how the ordering of functional estimation and privacy perturbation affects the accuracy of private functional mean estimation. We introduce the Fit-Then-Protect (FTP) and Protect-Then-Fit (PTF) pipelines and establish a GDP guarantee for Gaussian perturbations of function-valued summaries in separable Banach spaces, and compare the accuracy of the two pipelines under a common privacy budget. Our results show that, under the stated conditions, FTP achieves strictly smaller estimation error than PTF for sufficiently large sample sizes, with the asymptotic error ratio governed by the basis dimension in finite-dimensional Hilbert spaces and by the number of sampling locations in RKHSs. These findings highlight the importance of exploiting functional structure before privatization.

The simulation results further demonstrate that the relative performance of FTP and PTF depends on the observation regime. In the noiseless and low-noise settings, the empirical results are consistent with the theoretical comparisons and show a substantial advantage for FTP. As observation noise increases, however, this advantage becomes smaller and can eventually disappear or reverse. Furthermore, the constraint radius $B'$ plays an important role in the trade-off: increasing $B'$ reduces the probability of an active constraint but simultaneously increases the magnitude of the privacy noise. Thus, statistical estimation, sensitivity calibration, and observation noise must be considered jointly when designing private functional procedures.

Several limitations of the current analysis suggest directions for future work. First, our direct theoretical comparison between FTP and PTF focuses primarily on the noiseless setting. Although we derive noisy-data bounds for FTP and investigate observation noise empirically, a general theoretical characterization of the FTP--PTF comparison with nonzero observation noise remains open. Such a result would need to account simultaneously for observation-noise covariance, privacy-induced covariance, and the nonlinear effects of constrained functional estimation. A more complete noisy-data theory could clarify how the relative efficiency of the two orderings changes across observation regimes and identify conditions under which their performance becomes comparable or one procedure may become preferable.

Second, the comparison could be extended to more general sampling designs, function spaces, and privacy mechanisms. Our analysis assumes a common observation grid, and our simulations focus on polynomial spaces and Gaussian-kernel RKHSs. Studying irregular sampling and broader classes of function spaces, such as Sobolev spaces and RKHSs with different eigenvalue decay properties, would help clarify how sampling design and functional smoothness affect the relative accuracy of the two pipelines. In addition, while the current PTF construction provides a convenient GDP guarantee through perturbation at the observation level, coordinatewise perturbation need not be optimal for the geometry of functional data. Characterizing the optimal GDP design for PTF, including structured perturbations that exploit basis, kernel, or covariance information, is therefore an important open problem. Such mechanisms could potentially improve the utility of PTF and provide a more geometry-aware comparison with FTP.

Finally, the FTP/PTF distinction extends naturally beyond functional mean estimation. Functional regression, functional principal component analysis, covariance-function estimation, and other functional data analysis procedures involve transformations of functional observations that may be performed either before or after privacy perturbation. Investigating how privacy-estimation ordering interacts with these procedures, particularly when the estimation operator is nonlinear, could lead to a broader theory of when statistical structure should be extracted before privatization and when privacy should instead be imposed at an earlier stage of the analysis.

Overall, our results highlight that differential privacy for functional data is not only a problem of calibrating an appropriate noise process, but also a problem of deciding where privacy should be imposed within the statistical pipeline. The comparison between FTP and PTF illustrates that this design choice can have a substantial effect on estimation accuracy, while also showing that its impact depends on the observation-noise regime and the geometry of the underlying function space. These findings motivate further development of privacy mechanisms that jointly account for functional structure, observation noise, sampling design, and the geometry of the statistical representation.


\acks{This work was supported in part by NSF grants SES-2150615 and SES-260910.  

}

\newpage

\appendix
\section{Limitation of Pure Differential Privacy}
In this section, we briefly review pure differential privacy (DP) and the Laplace mechanism, then use a constant-function example to illustrate how the privacy framework affects the relative performance of FTP and PTF.

\begin{definition}[$\epsilon$-Differential Privacy]
A randomized mechanism 
$M:\mathcal{X}^n\to \mathcal{Y}$ satisfies $\epsilon$-DP if for all adjacent data sets $X\sim X'$ and all measurable sets $A\subseteq \mathcal{Y}$,
\[
\Pr(M(X)\in A)
\le
e^\epsilon \Pr(M(X')\in A)
\]
\end{definition}

The privacy parameter $\epsilon$ controls the privacy loss, with smaller values providing stronger privacy guarantees and weaker data usability.
\begin{proposition}[Laplace Mechanism, \citealp{dwork2006calibrating}]
    Given a statistic $S(X)\in \mathbb{R}^d$, the Laplace mechanism releases
    \[
    M(X)=S(X)+L,\qquad
    L_i \overset{iid}{\sim} \mathrm{Lap}\left(\frac{\Delta_1(S)}{\epsilon}\right), i=1,\ldots,d.
    \]
    Then $M(X)$ satisfies $\epsilon$-DP.
\end{proposition}

In infinite-dimensional settings, extending Laplace-type mechanisms becomes difficult. In particular, some work shows that there does not exist a natural infinite-dimensional analogue of elliptical Laplace distributions that achieves pure $\epsilon$-DP \citep{reimherr2019elliptical}.

\begin{example}[A Motivating Constant-Function Example]
\label{ex: A motivating constant-function example}
To illustrate why the ordering of estimation and privacy may matter, consider the simplified setting in which there is no observation noise ($\sigma_0=0$) and each individual function is constant, $f_i(t)\equiv c_i,\ |c_i|\le B$. The target parameter is the sample mean $\bar c=n^{-1}\sum_{i=1}^n c_i$. 

Under the FTP strategy, each constant function is first estimated exactly, the estimates are averaged, and a single privacy mechanism is applied to the resulting scalar mean. Since the $L_1$-sensitivity of the sample mean is $2B/n$, the Laplace mechanism releases
\[
\widetilde c_{\mathrm{FTP}}
=
\bar c + \eta,
\qquad
\eta\sim \operatorname{Lap}
\!\left(
\frac{2B}{n\epsilon}
\right),
\]
which satisfies pure $\epsilon$-DP.

Under PTF, a natural implementation first computes the sample mean at each observation time, $\bar x_j=n^{-1}\sum_{i=1}^n x_{ij}$, privatizes each coordinate independently, and then estimates the constant function from the privatized observations. Allocating the privacy budget equally across the $T$ coordinates gives
\[
\widetilde x_j
=
\bar x_j
+
\eta_j,
\qquad
\eta_j
\overset{\text{i.i.d.}}{\sim}
\operatorname{Lap}
\!\left(
\frac{2BT}{n\epsilon}
\right),
\]
where the factor $T$ follows from the sequential composition theorem. Averaging the privatized coordinates yields a variance that is exactly $T$ times larger than that of the FTP estimator.

In contrast, under $\mu$-Gaussian differential privacy, the composition rule is quadratic rather than linear. If each coordinate is released under $\mu/\sqrt{T}$-GDP, then averaging the resulting Gaussian perturbations exactly offsets the $\sqrt{T}$ increase in the noise standard deviation, so both FTP and this coordinate-wise PTF construction achieve the same privacy-induced variance order. This toy example illustrates that the relative statistical performance of FTP and PTF depends not only on the estimation problem, but also on the privacy framework itself.
\end{example}

\section{Proofs and Technical Details}
In this section, we provide the proofs and additional technical details for our theoretical results in Section~\ref{sec: Theoretical Results}.
\begin{lemma}[Cameron-Martin Lemma, \citealp{mirshani2019formal}]
    \label{lemma: Cameron-Martin Lemma}
    Let $\mathbb{B}$ be a real separable Banach space, and let $Z\sim\mathcal N_{\mathbb B}(0,C)$ be a centered Gaussian random element with its Cameron--Martin space $\mathcal H_C$. Let $S:\mathcal X^n\to\mathbb B$ be a measurable deterministic summary compatible with $Z$. For $h\in\mathcal H_C$, let $T_h$ denote the associated measurable Gaussian linear functional. For $\sigma>0$, denote $Q$ as the probability measure over $\mathbb{B}$ induced by $\sigma Z$, and $\Set{P_X:X\in \mathcal{X}^n}$ the family of probability measures over $\mathbb{B}$ induced by $S(X)+ \sigma Z$, then every measure $P_X$ has a density over $\mathbb{B}$ with respect to $Q$, which is given by
    \[
    \frac{dP_X}{dQ}(y)
    =
    \exp\left\{
    -\frac{1}{2\sigma^2}
    \left(
    \|S(X)\|_{\mathcal H_C}^2
    -2T_{S(X)}(y)
    \right)
    \right\},
    \qquad Q\text{-almost everywhere}.
    \]
\end{lemma}

Here, $T_h$ is a measurable linear extension of $g\mapsto\langle h,g\rangle_{\mathcal H_C}$, defined on a full-Gaussian-measure linear subspace containing $\mathcal H_C$. The correspondence $h\mapsto T_h$ is linear, and $T_h(Z)\sim N(0,\|h\|_{\mathcal H_C}^2).$ These properties follow from the construction of the Gaussian linear functionals associated with the Cameron--Martin space in \citet[Section~2.2, Lemma~3.1]{mirshani2019formal}.

\begin{proof}\textbf{of Theorem~\ref{thm: GDP_functional}.}
Fix any two neighboring datasets $X\sim X'$. Let $P_X$ and $P_{X'}$ denote the laws induced by $S(X)+\sigma Z$ and $S(X')+\sigma Z$, respectively, and write $d=S(X')-S(X)\in\mathcal H_C$. If $d=0$, then $P_X=P_{X'}$ and $T(P_X,P_{X'})(\alpha)=1-\alpha\ge G_\mu(\alpha)$. Hence, suppose that $\|d\|_{\mathcal H_C}>0$.

By Lemma~\ref{lemma: Cameron-Martin Lemma}, the log-likelihood ratio is
\begin{align*}
L(y)
&:=\log\frac{dP_{X'}}{dP_X}(y)\\
&=\log \left(\dfrac{d P_{X'}}{dQ}(y) / \dfrac{d P_X}{dQ}(y)\right)\\
&=
-\frac{1}{2\sigma^2}
\left(
\|S(X')\|_{\mathcal H_C}^2
-\|S(X)\|_{\mathcal H_C}^2
-2T_d(y)
\right)\\
&=
\frac{1}{\sigma^2}T_d(y-S(X))
-\frac{\|d\|_{\mathcal H_C}^2}{2\sigma^2},
\end{align*}
almost surely under either distribution.

Define $U(y)=T_d(y-S(X))$, the Gaussian and translation properties of $T_d$ give
\[
U(Y)\sim
\begin{cases}
N(0,\sigma^2\|d\|_{\mathcal H_C}^2),
& Y\sim P_X,\\[2pt]
N(\|d\|_{\mathcal H_C}^2,
  \sigma^2\|d\|_{\mathcal H_C}^2),
& Y\sim P_{X'}.
\end{cases}
\]

By the Neyman--Pearson lemma, for $\alpha\in(0,1)$, a most powerful level-$\alpha$ test is given by
\[
\psi(y) = \mathbf{1}(L(y)>\eta_\alpha)
=\mathbf{1}(U(y)>c_{\alpha}).
\]

Thus $U(y)$ is the test statistic and the Type I error of $\psi$ is
\begin{align*}
    \alpha_\psi=\mathbb{E}_{P_X}[\psi(y)]
    =P_X(U(y)>c_\alpha)
    =1-\Phi\Bra{\dfrac{c_\alpha}{\sigma\|d\|_{\mathcal{H}_C}}}.
\end{align*}

When $\alpha_\psi\leq \alpha$, it follows that $c_\alpha\geq \sigma\|d\|_{\mathcal H_C}\Phi^{-1}(1-\alpha)$, and the type II error of $\psi$ is
\begin{align*}
    \beta_\psi
    =1-\mathbb{E}_{P_{X'}}[\psi(y)]
    =1-P_{X'}(U(y)>c_\alpha)
    =\Phi\Bra{\dfrac{c_\alpha-\|d\|_{\mathcal{H}_C}^2}{\sigma\|d\|_{\mathcal{H}_C}}}.
\end{align*}

Combine the above results and recall $\|d\|_{\mathcal H_C}/\sigma \le \Delta_{\mathcal H_C}(S)/\sigma\le \mu$, we get
\begin{align*}
T(P_X,P_{X'})(\alpha)
&=\inf_{\psi}\Set{\beta_\psi|\alpha_\psi\le\alpha}\\
&=
\Phi\left(
\frac{\sigma\|d\|_{\mathcal H_C}\Phi^{-1}(1-\alpha)-\|d\|_{\mathcal H_C}^2}
     {\sigma\|d\|_{\mathcal H_C}}
\right)\\
&=
\Phi\left(
\Phi^{-1}(1-\alpha)
-\frac{\|d\|_{\mathcal H_C}}{\sigma}
\right)\\
&\ge
\Phi\left(\Phi^{-1}(1-\alpha)-\mu\right)
=
G_\mu(\alpha),
\end{align*}

The cases where $\alpha=0,1$ follow by continuity. Since this holds for every neighboring data sets, $\widetilde S(X)=S(X)+\sigma Z$ satisfies $\mu$-GDP.
\end{proof}

\begin{proof}\textbf{of Corollary~\ref{col: finite-dimension functional GDP}.}
Take $\mathbb B=\mathcal H_m$. The covariance operator of $Z_m=\sum_{k=1}^m\xi_k\phi_k$ satisfies
\[
C_m\ell=\mathbb E[\ell(Z_m)Z_m]
=\sum_{j,k=1}^m\mathbb E[\xi_j\xi_k]\ell(\phi_j)\phi_k
=\sum_{k=1}^m\ell(\phi_k)\phi_k.
\]
For $g\in\mathcal H_m$, by the Riesz representation theorem, the functional $\ell_g(f)=\langle f,g\rangle_{\mathcal H_m}$ satisfies $C_m\ell_g=g$. Hence $C_m(\mathcal H_m^*)=\mathcal H_m$, and
\[
\langle C_m\ell_g,C_m\ell_v\rangle_{\mathcal H_{C_m}}
=\ell_g(C_m\ell_v)=\langle g,v\rangle_{\mathcal H_m}.
\]
Thus $\mathcal H_{C_m}=\mathcal H_m$ isometrically; the covariance corresponds to the identity under the Riesz identification with the dual.

Orthonormality gives
\[
\|S(X)-S(X')\|_{\mathcal H_{C_m}}^2
=\int_0^1\bigl[(\theta(X)-\theta(X'))^\top\phi(t)\bigr]^2\,dt
=\|\theta(X)-\theta(X')\|_2^2.
\]
Apply Theorem~\ref{thm: GDP_functional} with sensitivity $\Delta_{\mathcal{H}_m}$. Also $S(X)+\sigma Z_m=(\theta(X)+\sigma\xi)^\top\phi$, where $\xi\sim N(0,I_m)$. In Algorithm~\ref{alg: FTP} each fitted coefficient vector has norm at most $B'$, so $\Delta_{\mathcal{H}_m}\le2B'/n$, proving the stated scale and implementation.
\end{proof}

\begin{proof}\textbf{of Corollary~\ref{col: RKHS GDP}.}
Realize the continuous version of $Z_K$ as a Borel Gaussian random element in $\mathbb B=C[0,1]$. Borel measurability follows because continuous paths are determined by evaluations on a countable dense set; continuous linear functionals can be approximated by finite linear combinations of evaluations. Every $\ell\in\mathbb B^*$ is represented by a finite signed Borel measure $\nu_\ell$ through $\ell(f)=\int_0^1 f(s)\,\nu_\ell(ds)$. Hence
\[
(C_K\ell)(t)=\mathbb E[\ell(Z_K)Z_K(t)]
=\int_0^1 K(t,s)\,\nu_\ell(ds).
\]
Fubini's theorem applies because $K$ is bounded on the compact domain and $\nu_\ell$ has finite total variation.

The kernel sections depend continuously on their location in RKHS norm, since
\[
\|K(s,\cdot)-K(t,\cdot)\|_{\mathcal H_K}^2
=K(s,s)+K(t,t)-2K(s,t).
\]
Thus $C_K\ell=\int K(s,\cdot)\,\nu_\ell(ds)$ is a well-defined $\mathcal H_K$-valued integral. For $\ell,g\in\mathbb B^*$,
\[
\langle C_K\ell,C_Kg\rangle_{\mathcal H_K}
=\int\!\int K(s,t)\,\nu_\ell(ds)\nu_g(dt)
=\ell(C_Kg)
=\langle C_K\ell,C_Kg\rangle_{\mathcal H_{C_K}}.
\]

The evaluation functional $\delta_t(f)=f(t)$ satisfies $C_K\delta_t=K(t,\cdot)$. Consequently $C_K(\mathbb B^*)\subseteq\mathcal H_K$ contains a dense span of kernel sections. Completing this range proves $\mathcal H_{C_K}=\mathcal H_K$ isometrically. The identification agrees inside $C[0,1]$, since $\|f\|_\infty\le\sup_{t\in[0,1]}\sqrt{K(t,t)}\,\|f\|_{\mathcal H_K}$.

Theorem~\ref{thm: GDP_functional} applies with sensitivity $\Delta_{\mathcal{H}_K}$, and Algorithm~\ref{alg: FTP} has $\Delta_{\mathcal{H}_K}\le2B'/n$. Grid evaluation is measurable post-processing, and the GP definition gives $(Z_K(t_1),\ldots,Z_K(t_T))^\top\sim N(0,K_T)$. This proves both the full-function guarantee and the grid implementation.
\end{proof}

\begin{proof}\textbf{of Proposition~\ref{prop: GDP for PTF}.}
Since $\widetilde{x}$ satisfies $\mu$-GDP by composition, and $\mathcal F_{\mathrm{PTF}}:\mathcal X\rightarrow\mathbb{B}$ is a deterministic measurable mapping, the result follows immediately from the post-processing property.
\end{proof}
\begin{proof}\textbf{of Theorem~\ref{thm: FTP IMSE}.}
Denote $d=\overline{\hat{\theta}^{(c)}}-\bar{\theta}\in \R^m$. Since $\Phi$ has full column rank, $M=\Phi^\top \Phi$ is positive
definite. Notice that the true coefficient vectors satisfy
$\|\theta_i\|_2 \leq B, \ i=1,\ldots,n$. Define the constraint set
\[
\mathcal C_{B'}
:=
\left\{
\theta\in\mathbb R^m:\|\theta\|_2\leq B'
\right\}.
\]
For each individual $i$, first consider the unconstrained least-squares
estimator
\[
\hat{\theta}_i^{(u)}
:=
\arg\min_{\theta\in\mathbb R^m}
\|x_i-\Phi\theta\|_2^2.
\]
We denote the following notations for convenience:
\begin{equation*}
    \begin{aligned}
        e_i = \hat{\theta}_i^{(c)} - \theta_i, &&
        \xi_i = \hat{\theta}_i^{(u)} - \theta_i, &&
        r_i = \hat{\theta}_i^{(c)} - \hat{\theta}_i^{(u)}, &&
        b_i = \E e_i.
    \end{aligned}
\end{equation*}
It follows that $e_i = \xi_i + r_i$. Since $M=\Phi^\top\Phi$ is positive definite, $\hat{\theta}_i^{(u)} = M^{-1}\Phi^\top x_i$. Using $x_i=\Phi\theta_i+\epsilon_i$, we obtain $\xi_i = \hat{\theta}_i^{(u)}-\theta_i = M^{-1}\Phi^\top\epsilon_i$. Because $\epsilon_i\sim N(0,\sigma_0^2I_T)$, it follows that $\xi_i\sim N(0,\sigma_0^2M^{-1})$. Consequently,
\begin{align}
\mathbb E
\left\|
\xi_i
\right\|_2^2
=
\operatorname{tr}
\left[
\operatorname{Cov}
\left(
\xi_i
\right)
\right]
=
\sigma_0^2\operatorname{tr}(M^{-1}).
\label{eq:unconstrained_mse}
\end{align}
We next express the constrained estimator as a metric projection.
\[
\hat{\theta}_i^{(c)}
=
\arg\min_{\theta\in\mathcal C_{B'}}
\|x_i-\Phi\theta\|_2^2.
\]
Using $\Phi^\top x_i=M\hat{\theta}_i^{(u)}$, we can complete the square:
\begin{align*}
\|x_i-\Phi\theta\|_2^2
&=
x_i^\top x_i
-2\theta^\top\Phi^\top x_i
+\theta^\top M\theta\\
&=
x_i^\top x_i
-2\theta^\top M\hat{\theta}_i^{(u)}
+\theta^\top M\theta\\
&=
C_i+
(\theta-\hat{\theta}_i^{(u)})^\top
M
(\theta-\hat{\theta}_i^{(u)}),
\end{align*}
where 
$C_i = x_i^\top x_i - \hat{\theta}_i^{(u)\top} M\hat{\theta}_i^{(u)}$ 
does not depend on $\theta$. Hence
\begin{equation}
\hat{\theta}_i^{(c)}
=
\arg\min_{\theta\in\mathcal C_{B'}}
(\theta-\hat{\theta}_i^{(u)})^\top
M
(\theta-\hat{\theta}_i^{(u)}).
\label{eq:constrained_projection}
\end{equation}
Define an inner product on $\mathbb R^m$ by 
$\langle u,v\rangle_M:=u^\top M v,\ u,v\in\mathbb R^m$,
and let 
$\|u\|_M:=\sqrt{\langle u,u\rangle_M}=\sqrt{u^\top Mu}$ denote the corresponding norm. Then $(\mathbb R^m,\langle\cdot,\cdot\rangle_M)$ is a Hilbert space since $M$ is positive definite. Equation~\eqref{eq:constrained_projection} can therefore be written as
\begin{equation}
\hat{\theta}_i^{(c)}
=\argmin_{\theta\in \mathcal{C}_{B'}} 
\left\|
\theta-\hat{\theta}_i^{(u)}
\right\|_{M}^2
=
P_{\mathcal C_{B'}}^{M}
\left(
\hat{\theta}_i^{(u)}
\right),
\label{eq:metric_projection}
\end{equation}
where $P_{\mathcal C_{B'}}^{M}$ denotes the metric projection onto
$\mathcal C_{B'}$ with respect to the $M$-norm.\vspace{1em}\\
The set $\mathcal C_{B'}$ is nonempty, closed, and convex. By
Proposition~4.16 of \citet{bauschke2017convexity},
the metric projection onto a nonempty closed convex subset of a
Hilbert space is firmly nonexpansive, and hence nonexpansive. Thus,
in the Hilbert space
$(\mathbb R^m,\langle\cdot,\cdot\rangle_M)$, for all $u,v\in\mathbb R^m$, we have
\begin{equation}
\left\|
P_{\mathcal C_{B'}}^{M}(u)
-
P_{\mathcal C_{B'}}^{M}(v)
\right\|_M
\leq
\|u-v\|_M.
\label{eq:projection_nonexpansive}
\end{equation}
Since $\|\theta_i\|_2\leq B\leq B'$, we have
$\theta_i\in\mathcal C_{B'}$. Therefore, $P_{\mathcal C_{B'}}^{M}(\theta_i)=\theta_i$. Applying~\eqref{eq:projection_nonexpansive} with $u=\hat{\theta}_i^{(u)},\ v=\theta_i$, and using~\eqref{eq:metric_projection}, we obtain
\begin{align}
\left\|
e_i\right\|_M
=
\left\|
\hat{\theta}_i^{(c)}-\theta_i
\right\|_M
=
\left\|
P_{\mathcal C_{B'}}^{M}
\left(
\hat{\theta}_i^{(u)}
\right)
-
P_{\mathcal C_{B'}}^{M}(\theta_i)
\right\|_M
\leq
\left\|
\hat{\theta}_i^{(u)}-\theta_i
\right\|_M
=
\left\|
\xi_i\right\|_M.
\label{eq:M_norm_nonexpansive}
\end{align}
Since $M$ is positive definite, for every $v\in\mathbb R^m$,
\begin{equation}
\lambda_{\min}(M)\|v\|_2^2
\leq
\|v\|_M^2
\leq
\lambda_{\max}(M)\|v\|_2^2.
\label{eq: M_norm_to_l2}
\end{equation}
Consequently,~\eqref{eq:M_norm_nonexpansive} implies
\begin{align*}
\lambda_{\min}(M)
\left\|
e_i
\right\|_2^2
\leq
\left\|
e_i
\right\|_M^2
\leq
\left\|
\xi_i
\right\|_M^2
\leq
\lambda_{\max}(M)
\left\|
\xi_i
\right\|_2^2.
\end{align*}
Therefore,
\begin{equation}
\left\|
e_i
\right\|_2^2
\leq
\kappa(M)
\left\|
\xi_i
\right\|_2^2,
\label{eq:condition_number_bound}
\end{equation}
where $\kappa(M)=\lambda_{\max}(M)/\lambda_{\min}(M)$ is the spectral condition number of $M$.\\
Taking expectations in~\eqref{eq:condition_number_bound} and applying
\eqref{eq:unconstrained_mse}, we obtain
\begin{equation}
\mathbb E
\left\|
e_i
\right\|_2^2
\leq
\kappa(M)\sigma_0^2
\operatorname{tr}(M^{-1}).
\label{eq:individual_constrained_bound}
\end{equation}
Denote $C = \kappa(M)\sigma_0^2
\operatorname{tr}(M^{-1})$, notice that $e_i$'s are independent
\begin{align}
\mathbb E\|d\|_2^2
=
\dfrac{1}{n^2} \mathbb E \left\|\sum_{i=1}^n e_i\right\|_2^2
=
\dfrac{1}{n^2} \mathbb E \Bra{
\sum_{i=1}^n \|e_i\|_2^2 + \sum_{i\ne j} e_i^\top e_j}
=
\dfrac{1}{n^2} \Bra{
\sum_{i=1}^n \mathbb E \|e_i\|_2^2
+ \sum_{i\ne j} \mathbb E e_i^\top \mathbb E e_j
}.
\label{eq: direct expand d}
\end{align}
Recall Equation~\eqref{eq:individual_constrained_bound}, it follows that
\begin{equation}
\dfrac{1}{n^2} \sum_{i=1}^n \mathbb E \|e_i\|_2^2 
\le 
\dfrac{\kappa(M)\sigma_0^2
\operatorname{tr}(M^{-1})}{n}.
\label{eq: noise term}
\end{equation}
Then, we only need to bound the clipping bias term.\\
Recall the definition of $\hat{\theta}_i^{(c)}$ given by Equation~\eqref{eq:metric_projection}, we have
\begin{equation*}
\|r_i\|_M
=
\left\|
\hat{\theta}_i^{(c)} - \hat{\theta}_i^{(u)}
\right\|_M
\le 
\left\|
\theta_i - \hat{\theta}_i^{(u)}
\right\|_M
=
\|\xi_i\|_M.
\end{equation*}
Combine with Equation~\eqref{eq: M_norm_to_l2}, we have $\|\hat{\theta}_i^{(c)} - \hat{\theta}_i^{(u)}\|_2^2 
\le 
\kappa(M) \|\hat{\theta}_i^{(u)} - \theta_i\|_2^2$. Since $\hat{\theta}_i^{(u)}$ is the unconstrained minimizer, we know that $r_i=0$ if $\|\hat{\theta}_i^{(u)}\|_2\le B'$. Therefore, we can rewrite 
$r_i = \Bra{
\hat{\theta}_i^{(c)} - \hat{\theta}_i^{(u)}
}\cdot \mathbf{1} \Set{\|\hat{\theta}_i^{(u)}\|_2>B'}$.
It follows that
\begin{align*}
\|r_i\|_2^2
=\|\hat{\theta}_i^{(c)} - \hat{\theta}_i^{(u)}\|_2^2 \cdot \mathbf{1} \Set{\|\hat{\theta}_i^{(u)}\|_2>B'}
\le
\kappa(M) \|\hat{\theta}_i^{(u)} - \theta_i\|_2^2
\cdot \mathbf{1} \Set{\|\hat{\theta}_i^{(u)}\|_2>B'}.
\end{align*}
Thus, we can obtain
\begin{align}
\|b_i\|_2^2
=
\| \mathbb E e_i \|_2^2
=
\| \mathbb E r_i \|_2^2
\le
\mathbb E \|r_i\|_2^2
\le
\kappa(M) \mathbb E \SqBra{
\|\hat{\theta}_i^{(u)} - \theta_i\|_2^2
\cdot \mathbf{1} \Set{\|\hat{\theta}_i^{(u)}\|_2>B'}
}.
\label{eq: clipping bias bound}
\end{align}
By the Cauchy--Schwarz inequality, we have
\begin{align}
\mathbb E \SqBra{
\|\hat{\theta}_i^{(u)} - \theta_i\|_2^2
\cdot \mathbf{1} \Set{\|\hat{\theta}_i^{(u)}\|_2>B'}
}
=
\mathbb E \SqBra{
\|\xi_i\|_2^2
\cdot \mathbf{1} \Set{\|\hat{\theta}_i^{(u)}\|_2>B'}
}
\nonumber\\
\le
\sqrt{\mathbb E \|\xi_i\|_2^4} \sqrt{\operatorname{Pr}(\|\hat{\theta}_i^{(u)}\|_2>B')}.
\label{eq: Cauchy--Schwarz inequality}
\end{align}
Since $\xi_i\overset{i.i.d.}{\sim} N(0,\Sigma)$ with $\Sigma = \sigma_0^2 M^{-1}$, we have
\begin{equation}
\sqrt{\mathbb E \|\xi_i\|_2^4}
=
\sqrt{(\operatorname{tr}\Sigma)^2 + 2\operatorname{tr}(\Sigma^{2})}
\le \sqrt{3 (\operatorname{tr}\Sigma)^2}
= \sqrt{3} \sigma_0^2 \operatorname{tr}(M^{-1}).
\label{eq: 4th moment guassian bound}
\end{equation}
Recall $\|\theta_i\|_2\le B$, denote $\delta=B'-B$, then 
$\Set{\|\hat{\theta}_i^{(u)}\|_2>B'}$ implies 
$\Set{\|\hat{\theta}_i^{(u)}-\theta_i\|_2>\delta}$, which leads to
\begin{align*}
\operatorname{Pr}(\|\hat{\theta}_i^{(u)}\|_2>B')
\le
\operatorname{Pr}(\|\hat{\theta}_i^{(u)}-\theta_i\|_2>\delta)
=\operatorname{Pr}(\|\xi_i\|_2>\delta).
\end{align*}
Let $\xi_i = \sigma_0 M^{-1/2}Z_i$, $Z_i\sim N(0, I_m)$, $\|Z_i\|_2^2\sim \chi^2_m$, we have
\[
\|\xi_i\|_2^2
= 
\sigma_0^2 Z_i^\top M^{-1} Z_i 
\le 
\sigma_0^2 \lambda_{\max}(M^{-1}) \|Z_i\|_2^2
=
\dfrac{\sigma_0^2}{\lambda_{\min}(M)} \|Z_i\|_2^2.
\]
Then, we can obtain
\[
\operatorname{Pr}(\|\xi_i\|_2>\delta)
=
\operatorname{Pr}(\|\xi_i\|_2^2>\delta^2)
\le 
\operatorname{Pr}\Bra{
\dfrac{\sigma_0^2}{\lambda_{\min}(M)} \|Z_i\|_2^2>\delta^2
}
=\operatorname{Pr}
\Bra{
\chi^2_m > \dfrac{\delta^2 \lambda_{\min}(M)}{\sigma_0^2}
}.
\]
By Lemma 1 from \citet{laurent2000adaptive},
\begin{equation}
\operatorname{Pr}(\chi^2_m > (\sqrt{m}+\sqrt{2x})^2)
\le
\operatorname{Pr}(\chi^2_m > m + 2\sqrt{mx} + 2x)\le \exp(-x).
\label{eq: chi_square bound}
\end{equation}
By choosing $\sqrt{m}+\sqrt{2x} = \dfrac{\delta\sqrt{\lambda_{\min}(M)}}{\sigma_0}$ (when $\delta>\dfrac{\sqrt{m}\sigma_0}{\sqrt{\lambda_{\min}(M)}}$), we have
\begin{align}
\operatorname{Pr}(\|\hat{\theta}_i^{(u)}\|_2>B')
&\le
\operatorname{Pr}(\|\xi_i\|_2>\delta)
\nonumber\\
&\le
\operatorname{Pr}
\Bra{
\chi^2_m > \dfrac{\delta^2 \lambda_{\min}(M)}{\sigma_0^2}
}
\nonumber\\
&\le
\exp\SqBra{
-\dfrac{1}{2}\Bra{
\dfrac{\delta\sqrt{\lambda_{\min}(M)}}{\sigma_0}
- \sqrt{m}
}^2
}.
\label{eq: clipping probability bound}
\end{align}
Combine Equation~\eqref{eq: clipping bias bound},~\eqref{eq: Cauchy--Schwarz inequality},~\eqref{eq: 4th moment guassian bound} and~\eqref{eq: clipping probability bound}, we obtain
\begin{align}
\|b_i\|_2^2
&\le
\kappa(M) \mathbb E \SqBra{
\|\hat{\theta}_i^{(u)} - \theta_i\|_2^2
\cdot \mathbf{1} \Set{\|\hat{\theta}_i^{(u)}\|_2>B'}
}
\nonumber\\
&\le
\kappa(M)
\sqrt{\mathbb E \|\xi_i\|_2^4} \sqrt{\operatorname{Pr}(\|\hat{\theta}_i^{(u)}\|_2>B')}
\nonumber\\
&\le
\kappa(M)
\sqrt{3} \sigma_0^2 \operatorname{tr}(M^{-1})
\exp\SqBra{
-\dfrac{1}{4}\Bra{
\dfrac{\delta\sqrt{\lambda_{\min}(M)}}{\sigma_0}
- \sqrt{m}
}^2
}.
\label{eq: individual_bias_term}
\end{align}
By the Cauchy--Schwarz inequality, we have
\begin{equation}
\dfrac{1}{n^2} \sum_{i\ne j} \mathbb E e_i^\top \mathbb E e_j
=
\dfrac{1}{n^2} \sum_{i\ne j} b_i^\top b_j
\le
\Bra{1-\dfrac{1}{n}} \max_{i} \|b_i\|^2_2.
\label{eq: Cauchy--Schwarz}
\end{equation}
Combine Equation~\eqref{eq: direct expand d},~\eqref{eq: noise term},~\eqref{eq: individual_bias_term} and~\eqref{eq: Cauchy--Schwarz}, we get
\begin{equation}
\mathbb E\|d\|_2^2
\leq
\dfrac{\kappa(M)\sigma_0^2
\operatorname{tr}(M^{-1})}{n}
+
\Bra{1-\dfrac{1}{n}}
\kappa(M)
\sqrt{3} \sigma_0^2 \operatorname{tr}(M^{-1})
\exp\SqBra{
-\dfrac{1}{4}\Bra{
\dfrac{\delta\sqrt{\lambda_{\min}(M)}}{\sigma_0}
- \sqrt{m}
}^2
}.
\label{eq:d_final_bound}
\end{equation}
Finally, since $\widetilde{\theta}_{\mathrm{FTP}}
=
\overline{\hat{\theta}^{(c)}}+Z
=
\bar{\theta}+d+Z$, and $\int_0^1\phi(t)\phi(t)^\top\,dt=I_m$, we have
\begin{align*}
\mathrm{IMSE}_{\mathrm{FTP}}
=
\mathbb E
\int_0^1
\left[
{\widetilde{\theta}_{\mathrm{FTP}}}^\top\phi(t)
-
\bar{\theta}^\top\phi(t)
\right]^2dt
=
\mathbb E\|d+Z\|_2^2.
\end{align*}
Because $Z\sim N(0,\sigma^2I_m)$ is independent of $d$ and
$\mathbb E[Z]=0$,
\[
\mathbb E\|d+Z\|_2^2
=
\mathbb E\|d\|_2^2
+
\mathbb E\|Z\|_2^2
=
\mathbb E\|d\|_2^2+m\sigma^2.
\]
With $\sigma=2B'/(n\mu)$, we have $m\sigma^2=4m{B'}^2/(n^2\mu^2)$.\\
Combining this identity with~\eqref{eq:d_final_bound} yields
\begin{equation}
\begin{aligned}
\mathrm{IMSE}_{\mathrm{FTP}}
\leq
&\frac{4m{B'}^2}{n^2\mu^2}
+
\dfrac{\kappa(M)\sigma_0^2
\operatorname{tr}(M^{-1})}{n}\\
&+
\Bra{1-\dfrac{1}{n}}
\kappa(M)
\sqrt{3} \sigma_0^2 \operatorname{tr}(M^{-1})
\exp\SqBra{
-\dfrac{1}{4}\Bra{
\dfrac{(B'-B)\sqrt{\lambda_{\min}(M)}}{\sigma_0}
- \sqrt{m}
}^2
}.
\end{aligned}
\label{eq: final_IMSE_FTP bound}
\end{equation}
When $\sigma_0=0$ and $B'\ge B$, the true coefficient vector is feasible and achieves zero fitting loss. Since $\Phi$ has full column rank, $\hat\theta_i^{(c)}=\theta_i$. Thus the estimation error consists only of the privacy perturbation, giving $\operatorname{IMSE}_{\mathrm{FTP}}=4m(B')^2/(n^2\mu^2)$.
\end{proof}
\begin{proof}\textbf{of Theorem~\ref{thm: PTF IMSE lower}.}
Denote $a_j=\|\phi(t_j)\|_2$, when $\sigma_0=0$, $x_i=\Phi\theta_i$.
Since $\| \theta_i\|_2\le B<B', i=1,\ldots,n$, we have
\[
|x_{ij}|
=
|\phi(t_j)^\top\theta_i|
\leq
Ba_j
<
B'a_j:=c_j.
\]
Therefore the clipping $\operatorname{clip}(x_{ij},c_j)$ is inactive and
\[
\bar x^{(\mathrm{clip})}
=
\bar{x}
=
\dfrac{1}{n}\sum_{i=1}^n x_i
=
\Phi\bar\theta.
\]
It follows that the released private vector can be written as $\widetilde x
=
\Phi\bar\theta+Z$, where
\[
Z\sim N(0,\tau^2W^{-1}),
\quad
\tau^2=\frac{4B'^2T}{n^2\mu^2},
\quad
W = \mathrm{diag}\{a_1^{-2},\ldots,a_T^{-2}\}.
\]
We first analyze the corresponding unconstrained WLS estimator. By
definition,
\[
\widetilde\theta_{\mathrm{PTF}}^{(u)}
=
A^{-1}\Phi^\top W\widetilde x.
\]
Substituting $\widetilde x=\Phi\bar\theta+Z$ gives
\[
\widetilde\theta_{\mathrm{PTF}}^{(u)}
=
\bar\theta+
A^{-1}\Phi^\top WZ.
\]
Since $\operatorname{Cov}(Z)=\tau^2W^{-1}$, we obtain
\[
\begin{aligned}
\operatorname{Cov}
\left(
\widetilde\theta_{\mathrm{PTF}}^{(u)}
\right)
&=
A^{-1}\Phi^\top W
(\tau^2W^{-1})
W\Phi A^{-1}\\
&=
\tau^2
A^{-1}\Phi^\top W\Phi A^{-1}\\
&=
\tau^2A^{-1}.
\end{aligned}
\]
Hence $\widetilde\theta_{\mathrm{PTF}}^{(u)}-\bar\theta
\sim
N(0,\tau^2A^{-1})$, and therefore
\begin{equation}
\mathbb E
\left\|
\widetilde\theta_{\mathrm{PTF}}^{(u)}
-\bar\theta
\right\|_2^2
=
\tau^2\operatorname{tr}(A^{-1}).
\label{eq: PTF unconstrained MSE}
\end{equation}
We next account for the constraint $\|\theta\|_2\le B'$. Let $\mathcal A=\{\|\widetilde\theta_{\mathrm{PTF}}^{(u)}\|_2>B'\}$ denote the event that the unconstrained estimator violates the
constraint. On the complementary event $\mathcal A^c$, the
unconstrained estimator is feasible and hence
$\widetilde\theta_{\mathrm{PTF}}^{(c)}
=
\widetilde\theta_{\mathrm{PTF}}^{(u)}$.
Therefore
\[
\begin{aligned}
\mathbb E\|\widetilde\theta_{\mathrm{PTF}}^{(c)}-\bar\theta\|_2^2
&=
\mathbb E[
\|\widetilde\theta_{\mathrm{PTF}}^{(u)}-\bar\theta\|_2^2
\mathbf 1_{\mathcal A^c}]
+
\mathbb E[
\|\widetilde\theta_{\mathrm{PTF}}^{(c)}-\bar\theta\|_2^2
\mathbf 1_{\mathcal A}]
\\
&\geq
\mathbb E[
\|\widetilde\theta_{\mathrm{PTF}}^{(u)}-\bar\theta\|_2^2
\mathbf 1_{\mathcal A^c}]
\\
&=
\mathbb E\|\widetilde\theta_{\mathrm{PTF}}^{(u)}-\bar\theta\|_2^2
-
\mathbb E[
\|\widetilde\theta_{\mathrm{PTF}}^{(u)}-\bar\theta\|_2^2
\mathbf 1_{\mathcal A}].
\end{aligned}
\]
Thus, by Equation~\eqref{eq: PTF unconstrained MSE}, we get
\begin{equation}
\operatorname{IMSE}_{\mathrm{PTF}}
\geq
\tau^2\operatorname{tr}(A^{-1})
-
R_{\mathrm{PTF}},
\label{eq: PTF remainder decomposition}
\end{equation}
where $R_{\mathrm{PTF}}
=
\mathbb E[
\|\widetilde\theta_{\mathrm{PTF}}^{(u)}-\bar\theta\|_2^2
\mathbf 1_{\mathcal A}]$. It remains to bound the remainder term. Let $\xi=\widetilde\theta_{\mathrm{PTF}}^{(u)}-\bar\theta$. Then $\xi\sim N(0,\tau^2A^{-1})$. Since $\|\bar\theta\|_2\le B$, the event $\mathcal A$ implies
\[
\|\xi\|_2
=
\|\widetilde\theta_{\mathrm{PTF}}^{(u)}-\bar\theta\|_2
\geq
\|\widetilde\theta_{\mathrm{PTF}}^{(u)}\|_2-\|\bar\theta\|_2
>
B'-B.
\]
Consequently,
\[
\mathcal A
\subseteq
\left\{
\|\xi\|_2>B'-B
\right\}.
\]
Write $\xi=\tau A^{-1/2}z,
\ 
z\sim N(0,I_m)$. Since
\[
\|A^{-1/2}z\|_2^2
\leq
\lambda_{\max}(A^{-1})\|z\|_2^2
=
\frac{\|z\|_2^2}
{\lambda_{\min}(A)},
\]
we obtain
\[
\Pr(\mathcal A)
\leq
\Pr\left(
\|z\|_2^2
>
\frac{(B'-B)^2\lambda_{\min}(A)}
{\tau^2}
\right).
\]
Because $\|z\|_2^2\sim\chi_m^2$, apply the chi-square tail bound in Equation~\eqref{eq: chi_square bound}, we get
\begin{equation}
\Pr(\mathcal A)
\leq
\exp\left[
-\frac12
\left\{
\frac{(B'-B)\sqrt{\lambda_{\min}(A)}}{\tau}
-\sqrt m
\right\}^2
\right].
\label{eq: PTF active probability}
\end{equation}
Finally, by the Cauchy--Schwarz inequality,
\[
R_{\mathrm{PTF}}
\leq
\left(
\mathbb E\|\xi\|_2^4
\right)^{1/2}
\Pr(\mathcal A)^{1/2}.
\]
For a centered Gaussian vector with covariance
$\tau^2A^{-1}$,
\[
\mathbb E\|\xi\|_2^4
=
\tau^4
\left[
\{\operatorname{tr}(A^{-1})\}^2
+
2\operatorname{tr}(A^{-2})
\right]
\le 
\tau^4 [3\{\operatorname{tr}(A^{-1})\}^2].
\]
Therefore,
\[
\left(
\mathbb E\|\xi\|_2^4
\right)^{1/2}
\leq
\tau^2
\sqrt{3}\operatorname{tr}(A^{-1}).
\]
Combining this inequality with Equation~\eqref{eq: PTF active probability}
and Equation~\eqref{eq: PTF remainder decomposition} gives
Equation~\eqref{eq: PTF lower bound}.
\end{proof}

\begin{proof}
\textbf{of Corollary~\ref{col: IMSE Ratio for Polynomials}.}
Write $\tau_n=\tau=2B_n'\sqrt T/(n\mu)$. Since $B_n'$ is a non-decreasing sequence  satisfying $B_n'>B$ for sufficiently large $n$, we have 
\[
\frac{B_n'-B}{\tau_n}
=
\frac{n\mu}{2\sqrt T}
\left(1-\frac{B}{B_n'}\right)
\longrightarrow\infty.
\]
Thus, the radius condition in Theorem~\ref{thm: PTF IMSE lower} holds for all sufficiently large $n$, and its exponential remainder tends to zero. Combining that lower bound with the exact noiseless FTP error in Theorem~\ref{thm: FTP IMSE} gives
\[
\limsup_{n\to\infty}
\frac{\operatorname{IMSE}(\widetilde f_{\mathrm{FTP}})}
     {\operatorname{IMSE}(\widetilde f_{\mathrm{PTF}})}
\le
\frac{m}{T\operatorname{tr}(A^{-1})}.
\]

Define $u_j=\phi(t_j)/\|\phi(t_j)\|_2\in\mathbb R^m$ and $A_j=u_ju_j^\top$. Then $A_j$ is a rank-one orthogonal projector with $\operatorname{tr}(A_j)=1$. By the definition of $W$, $A=\Phi^\top W\Phi=\sum_{j=1}^T A_j$, so $\operatorname{tr}(A)=T$. Let $\lambda_1,\ldots,\lambda_m>0$ be the eigenvalues of $A$. Then
\[
\frac{m}{T\operatorname{tr}(A^{-1})}
=
\frac{m}{T\sum_{k=1}^m\lambda_k^{-1}}
\le
\frac{\sum_{k=1}^m\lambda_k}{Tm}
=
\frac1m.
\]
The inequality follows from the harmonic mean--arithmetic mean inequality, with equality if and only if all $\lambda_k$ are equal.
\end{proof}

\begin{proof}\textbf{of Theorem~\ref{thm: FTP MSE RKHS}.}
Let $\bar{\mathbf f}=n^{-1}\sum_{i=1}^n\mathbf f_i$. Since $K_T$ is positive definite, $K_T^{-1}$ and $K_T^{1/2}$ are well-defined.\\
\textcolor{red}{a.} 
We first show that for any $f_i\in \mathcal{H}_K$ with $\|f_i\|_{\mathcal{H}_K}\leq B$, we can always find another function $g_i$ in the span of $\{K(t_1,t),\ldots,K(t_T,t)\}$ such that 
\[
f_i(t_j)=g_i(t_j),\ \|g_i\|_{\mathcal{H}_K}\leq B.
\]
For each $i$, consider the minimum-norm interpolation problem
\[
g_i
=
\argmin_{g\in\mathcal H_K}
\left\{
\|g\|_{\mathcal H_K}:
g(t_j)=f_i(t_j),\ j=1,\ldots,T
\right\}.
\]
By the representer theorem,
\[
g_i(t)
=
\sum_{j=1}^T
(\gamma_i)_jK(t_j,t),
\quad
\gamma_i=K_T^{-1}\mathbf f_i.
\]
Hence
\[
\|g_i\|_{\mathcal H_K}^2
=
\gamma_i^\top K_T\gamma_i
=
\mathbf f_i^\top K_T^{-1}\mathbf f_i.
\]
Since $f_i$ itself is feasible for the interpolation problem,
\[
\|g_i\|_{\mathcal H_K}
\le
\|f_i\|_{\mathcal H_K}
\le B.
\]
Therefore,
\begin{equation}
\mathbf f_i^\top K_T^{-1}\mathbf f_i
\le B^2.
\label{eq: RKHS_min_norm_true_bound}
\end{equation}
Define
\[
\mathbf \beta_i
:=
K_T^{-1/2}\mathbf f_i.
\]
Then Equation~\eqref{eq: RKHS_min_norm_true_bound} gives
\begin{equation}
\|g_i\|_{\mathcal H_K} = \|\beta_i\|_2\le B.
\label{eq:true_transformed_norm}
\end{equation}
\textcolor{red}{b.} 
We then derive the form of the constrained and unconstrained solution to \begin{equation}
    \hat{f}_i^{(c)} 
    = \argmin_{\|f\|_{\mathcal{H}_K}\le B'} 
    \dfrac{1}{T}\sum_{j=1}^T [x_{ij}-f(t_j)]^2.
\label{eq: constrained_FTP_RKHS}
\end{equation}
and perform coefficient transformation for the convenience of following proof.\\ 
By the Representer theorem, the constrained estimator has the form
\[
\hat f_i^{(c)}(t)
=
\sum_{j=1}^T
(\hat\alpha_i^{(c)})_jK(t_j,t),
\]
and its RKHS norm satisfies
\[
\|\hat f_i^{(c)}\|_{\mathcal H_K}^2
=
\hat\alpha_i^{(c)\top}
K_T
\hat\alpha_i^{(c)}.
\]
Consequently, Equation~\eqref{eq: constrained_FTP_RKHS} is equivalent to
\begin{equation}
    \hat\alpha_i^{(c)}
    =
    \argmin_{\alpha^\top K_T\alpha\le {B'}^2}
    \|x_i-K_T\alpha\|_2^2.
    \label{eq: constrained_alpha_RKHS}
\end{equation}
Consider the corresponding unconstrained estimator
\[
\hat f_i^{(u)}
=
\argmin_{f\in\mathcal H_K}
\frac1T\sum_{j=1}^T
[x_{ij}-f(t_j)]^2.
\]
The minimum-norm unconstrained least-squares solution has the form
\begin{equation}
    \begin{aligned}
        \hat f_i^{(u)}(t)
        &=
        \sum_{j=1}^T [\hat\alpha_i^{(u)}]_jK(t_j,t),\\
        \hat\alpha_i^{(u)}
        &=
        K_T^{-1}x_i.
    \end{aligned}
\end{equation}
Indeed, the fitted-value vector is
\[
\hat{\mathbf{f}}_i^{(u)}
=
K_T\hat\alpha_i^{(u)}
=
x_i.
\]
Define
\[
\hat{\beta}_i^{(u)}
:=
K_T^{1/2}\hat\alpha_i^{(u)}
=
K_T^{-1/2}x_i
\text{\quad and\quad}
\hat{\beta} := K_T^{1/2}\alpha.
\]
We have
\[
\alpha^\top K_T\alpha\le {B'}^2
\Longleftrightarrow
\|\beta\|_2 = \|K_T^{1/2}\alpha\|_2\le B'.
\]
Hence, we can rewrite the optimization problem in Equation~\eqref{eq: constrained_alpha_RKHS} as
\begin{equation}
    \hat\beta_i^{(c)}
    =
    \argmin_{\|\beta\|_2\le B'}
    \|x_i-K_T^{1/2}\beta\|_2^2,
    \quad \hat\beta_i^{(c)} = K_T^{1/2}\hat\alpha_i^{(c)}.
    \label{eq: constrained_beta_RKHS}
\end{equation}
It follows
\[
\hat f_i^{(c)}(t)
=
\sum_{j=1}^T
(\hat\alpha_i^{(c)})_jK(t_j,t)
=
\sum_{j=1}^T
(K_T^{-1/2}\hat\beta_i^{(c)})_jK(t_j,t).
\]
\textcolor{red}{c.} 
Next, we introduce the metric projection and follow similar steps in the proof of Theorem~\ref{thm: FTP IMSE}.\\
Denote
\[
\mathcal C_{B'}
:=
\left\{
\beta\in\mathbb R^T: \|\beta\|_2\leq B'
\right\}.
\]
By substituting $\Phi = K_T^{1/2}$, $M = \Phi^\top \Phi = \Phi^2 = K_T$, we have
\begin{equation}
\hat{\beta}_i^{(c)}
=
\argmin_{\beta\in\mathcal C_{B'}}
(\beta-\hat{\beta}_i^{(u)})^\top
M
(\beta-\hat{\beta}_i^{(u)})
:= P_{\mathcal C_{B'}}^{M}
\left(
\hat{\beta}_i^{(u)}
\right).
\label{eq:constrained_projection_RKHS}
\end{equation}
where $P_{\mathcal C_{B'}}^{M}$ denotes the metric projection onto
$\mathcal C_{B'}$ with respect to the $M$-norm.\\
Since the projection is nonexpansive, and $\|\beta_i\|_2\leq B\leq B'$, we have
\begin{equation}
    \begin{aligned}
        \left\|
        \hat{\beta}_i^{(c)}-\beta_i
        \right\|_M
        =
        \left\|
        P_{\mathcal C_{B'}}^{M}
        \left(
        \hat{\beta}_i^{(u)}
        \right)
        -
        P_{\mathcal C_{B'}}^{M}(\beta_i)
        \right\|_M
        \leq
        \left\|
        \hat{\beta}_i^{(u)}-\beta_i
        \right\|_M.
    \end{aligned}
    \label{eq: M_norm_nonexpansive_RKHS}
\end{equation}
\textcolor{red}{d.} 
We denote convenient notations as follows
\begin{equation}
    \begin{aligned}
        e_i = \hat{\beta}_i^{(c)} - \beta_i, &&
        \xi_i = \hat{\beta}_i^{(u)} - \beta_i, &&
        r_i = \hat{\beta}_i^{(c)} - \hat{\beta}_i^{(u)}, &&
        b_i = \E e_i.
    \end{aligned}
    \label{eq: notations_RKHS}
\end{equation}
Then Equation~\eqref{eq: M_norm_nonexpansive_RKHS} becomes
\begin{equation}
    \|e_i\|_M\leq \|\xi_i\|_M.
    \label{eq: M_norm_nonexpansive_RKHS_simplified}
\end{equation}
Since
\[
x_i=\mathbf f_i+\epsilon_i,\quad
\epsilon_i\sim N(0,\sigma_0^2 I_T),
\]
we have
\[
\beta_i^{(u)}
=K_T^{-1/2} x_i
=
\Phi^{-1}x_i
=
\beta_i+\xi_i,
\]
where
\[
\xi_i
=
\Phi^{-1}\epsilon_i,\quad
\xi_i
\sim
N(0,\sigma_0^2M^{-1}).
\]
Consequently,
\begin{align}
\mathbb E
\left\|
\xi_i
\right\|_2^2
=
\operatorname{tr}
\left[
\operatorname{Cov}
\left(
\xi_i
\right)
\right]
=
\sigma_0^2\operatorname{tr}(M^{-1}),
\quad
\mathbb E
\left\|
\epsilon_i
\right\|_2^2
=\sigma_0^2 T.
\label{eq: unconstrained_mse_RKHS}
\end{align}
Then we can expand the MSE as
\begin{align}
\mathrm{MSE}_{\mathrm{FTP}}
&= 
\dfrac{1}{T} \sum_{j=1}^T \mathbb{E} 
[(\widetilde{f}_{\mathrm{FTP}}(t_j)-\bar{f}(t_j))^2]\nonumber\\
&= 
\dfrac{1}{T} \mathbb{E} 
\|\widetilde{\mathbf{f}}_{\mathrm{FTP}} - \bar{\mathbf{f}}\|_2^2\nonumber\\
&= 
\dfrac{1}{T} \mathbb{E} 
\|\overline{\hat{\mathbf{f}}^{(c)}} - \bar{\mathbf{f}} + \sigma \mathbf{Z}\|_2^2
\label{eq: direct expand MSE}.
\end{align}
Notice that $Z$ is a mean zero Gaussian process with kernel $K$, we have
\[
\mathbf{Z} = (Z(t_1),\ldots,Z(t_T))^\top
\sim \mathcal{N}(0,M),
\]
\[
\mathbb E\|\mathbf{Z}\|_2^2 = \operatorname{tr}(M).
\]
Thus, we have
\begin{align}
\mathrm{MSE}_{\mathrm{FTP}}
&= 
\dfrac{1}{T} \mathbb{E} 
\|\overline{\hat{\mathbf{f}}^{(c)}} - \bar{\mathbf{f}} + \sigma \mathbf{Z}\|_2^2\nonumber\\
&= 
\dfrac{1}{T} 
\Bra{\mathbb{E} 
\|\overline{\hat{\mathbf{f}}^{(c)}} - \bar{\mathbf{f}}\|_2^2 
+ 
\mathbb{E} \|\sigma \mathbf{Z}\|_2^2 + 0}
\nonumber\\
&= 
\dfrac{1}{T} 
\mathbb{E} 
\|\overline{\hat{\mathbf{f}}^{(c)}} - \bar{\mathbf{f}}\|_2^2 
+ 
\dfrac{4{B'}^2}{n^2\mu^2 T} \operatorname{tr}(M).
\label{eq: expanded MSE}
\end{align}
Note that $e_i$'s are independent, we can further expand the first term as
\begin{equation}
    \begin{aligned}
        \mathbb{E} 
        \|\overline{\hat{\mathbf{f}}^{(c)}} - \bar{\mathbf{f}}\|_2^2 
        &=
        \mathbb{E} 
        \|\overline{\hat{\mathbf{f}}^{(c)}} - \bar{\mathbf{g}}\|_2^2\\
        &=
        \mathbb{E} 
        \|\Phi \overline{\hat{\beta}^{(c)}} - \Phi\bar{\beta}\|_2^2\\
        &=
        \mathbb{E} 
        \|\overline{\hat{\beta}^{(c)}} - \bar{\beta}\|_{M}^2\\
        &=
        \dfrac{1}{n^2} \mathbb E \Bra{
        \sum_{i=1}^n \|e_i\|_M^2 + \sum_{i\ne j} e_i^\top  M e_j}\\
        &=
        \dfrac{1}{n^2} \Bra{
        \sum_{i=1}^n \mathbb E \|e_i\|_M^2
        + \sum_{i\ne j} \mathbb E e_i^\top M \mathbb E e_j
        }.
    \end{aligned}
    \label{eq: expand bias and clipping}
\end{equation}
Recall Equation~\eqref{eq: M_norm_nonexpansive_RKHS_simplified} and Equation~\eqref{eq: unconstrained_mse_RKHS}, we have
\begin{equation}
    \begin{aligned}
        \dfrac{1}{n^2} \sum_{i=1}^n \mathbb E \|e_i\|_M^2 
        \leq
        \dfrac{1}{n^2} \sum_{i=1}^n \mathbb E \|\xi_i\|_M^2
        =
        \dfrac{1}{n^2} \sum_{i=1}^n \mathbb E \|\epsilon_i\|_2^2
        =
        \dfrac{\sigma_0^2T}{n}.
    \end{aligned}
\label{eq: noise term_RKHS}
\end{equation}
\textcolor{red}{e.} 
Finally, we analyze the cross-term.\\
Recall the definition of $\hat{\beta}_i^{(c)}$ given by Equation~\eqref{eq:constrained_projection_RKHS}, we have
\begin{equation}
\|r_i\|_M
=
\left\|
\hat{\beta}_i^{(c)} - \hat{\beta}_i^{(u)}
\right\|_M
\le 
\left\|
\beta_i - \hat{\beta}_i^{(u)}
\right\|_M
=
\|\xi_i\|_M.
\end{equation}
Since $\hat{\beta}_i^{(u)}$ is the unconstrained minimizer, we know that
$r_i=0$ if $\|\hat{\beta}_i^{(u)}\|_2\le B'$.\\
Therefore, we can rewrite
\[
r_i 
= \Bra{
\hat{\beta}_i^{(c)} - \hat{\beta}_i^{(u)}
}\cdot \mathbf{1} \Set{\|\hat{\beta}_i^{(u)}\|_2>B'}.
\]
It follows that
\begin{align*}
\|r_i\|_M^2
=&\|\hat{\beta}_i^{(c)} - \hat{\beta}_i^{(u)}\|_M^2 \cdot \mathbf{1} \Set{\|\hat{\beta}_i^{(u)}\|_2>B'}\\
\le&
\|\xi_i\|_M^2
\cdot \mathbf{1} \Set{\|\hat{\beta}_i^{(u)}\|_2>B'}\\
=& \|\epsilon_i\|_2^2
\cdot \mathbf{1} \Set{\|\hat{\beta}_i^{(u)}\|_2>B'}.
\end{align*}
Thus, we can obtain that
\begin{equation}
    \|b_i\|_M^2
    =
    \| \mathbb E e_i \|_M^2
    =
    \| \mathbb E r_i \|_M^2
    \le
    \mathbb E \|r_i\|_M^2
    \le
    \mathbb E \SqBra{
    \|\epsilon_i\|_2^2
    \cdot \mathbf{1} \Set{\|\hat{\beta}_i^{(u)}\|_2>B'}
    }.
\label{eq: clipping bias bound_RKHS}
\end{equation}
By the Cauchy--Schwarz inequality, we have
\begin{equation}
    \mathbb E \SqBra{
    \|\epsilon_i\|_2^2
    \cdot \mathbf{1} \Set{\|\hat{\beta}_i^{(u)}\|_2>B'}
    }
    \leq
    \sqrt{\mathbb E \|\epsilon_i\|_2^4} \sqrt{\operatorname{Pr}(\|\hat{\beta}_i^{(u)}\|_2>B')}.
    \label{eq: Cauchy--Schwarz inequality_RKHS}
\end{equation}
Since $\epsilon_i\overset{i.i.d.}{\sim} N(0,\Sigma)$ with $\Sigma = \sigma_0^2 I_T$, we have
\begin{equation}
\sqrt{\mathbb E \|\epsilon_i\|_2^4}
=
\sqrt{(\operatorname{tr}\Sigma)^2 + 2\operatorname{tr}(\Sigma^{2})}
= \sigma_0^2 \sqrt{T^2+2T}.
\label{eq: 4th moment guassian bound_RKHS}
\end{equation}
Recall $\|\beta_i\|_2\le B$, denote $\delta=B'-B$, then 
$\Set{\|\hat{\beta}_i^{(u)}\|_2>B'}$ implies 
$\Set{\|\hat{\beta}_i^{(u)}-\beta_i\|_2>\delta}$, thus
\begin{align*}
\operatorname{Pr}(\|\hat{\beta}_i^{(u)}\|_2>B')
&\le
\operatorname{Pr}(\|\hat{\beta}_i^{(u)}-\beta_i\|_2>\delta)
\\
&=\operatorname{Pr}(\|\xi_i\|_2>\delta).
\end{align*}
Let $\xi_i = \sigma_0 M^{-1/2}Z_i$, $Z_i\sim N(0, I_T)$, $\|Z_i\|_2^2\sim \chi^2_T$, we have
\[
\|\xi_i\|_2^2
= 
\sigma_0^2 Z_i^\top M^{-1} Z_i 
\le 
\sigma_0^2 \lambda_{\max}(M^{-1}) \|Z_i\|_2^2
=
\dfrac{\sigma_0^2}{\lambda_{\min}(M)} \|Z_i\|_2^2.
\]
Apply the $\chi^2$ tail bound from \citet{laurent2000adaptive}, we have
\begin{align}
\operatorname{Pr}(\|\hat{\beta}_i^{(u)}\|_2>B')
&\le
\operatorname{Pr}(\|\xi_i\|_2>\delta)
\nonumber\\
&\le
\operatorname{Pr}
\Bra{
\chi^2_T > \dfrac{\delta^2 \lambda_{\min}(M)}{\sigma_0^2}
}
\nonumber\\
&\le
\exp\SqBra{
-\dfrac{1}{2}\Bra{
\dfrac{\delta\sqrt{\lambda_{\min}(M)}}{\sigma_0}
- \sqrt{T}
}^2
}.
\label{eq: clipping probability bound_RKHS}
\end{align}
Combine Equation~\eqref{eq: clipping bias bound_RKHS},~\eqref{eq: Cauchy--Schwarz inequality_RKHS},~\eqref{eq: 4th moment guassian bound_RKHS} and~\eqref{eq: clipping probability bound_RKHS}, we obtain
\begin{equation}
    \|b_i\|_M^2
    \leq
    \sigma_0^2\sqrt{T^2+2T}\cdot
    \exp\SqBra{
    -\dfrac{1}{4}\Bra{
    \dfrac{\delta\sqrt{\lambda_{\min}(M)}}{\sigma_0}
    - \sqrt{T}
    }^2
    }.
    \label{eq: individual_bias_term_RKHS}
\end{equation}
By the Cauchy--Schwarz inequality, we have
\begin{equation}
\dfrac{1}{n^2} \sum_{i\ne j} \mathbb E e_i^\top M \mathbb E e_j
=
\dfrac{1}{n^2} \sum_{i\ne j} b_i^\top M b_j
\le
\Bra{1-\dfrac{1}{n}} \max_{i} \|b_i\|^2_M.
\label{eq: Cauchy--Schwarz_RKHS}
\end{equation}
Finally combine Equation~\eqref{eq: expanded MSE},~\eqref{eq: expand bias and clipping},~\eqref{eq: noise term_RKHS},~\eqref{eq: individual_bias_term_RKHS},~\eqref{eq: Cauchy--Schwarz_RKHS}, and recall $K_T=M$, we get
\begin{equation}
\begin{aligned}
\mathrm{MSE}_{\mathrm{FTP}}
\leq
&\frac{4{B'}^2\operatorname{tr}(K_T)}{n^2\mu^2T}
+
\dfrac{\sigma_0^2}{n}\\
&+
\Bra{1-\dfrac{1}{n}}
\sigma_0^2 \sqrt{\dfrac{T+2}{T}}
\exp\SqBra{
-\dfrac{1}{4}\Bra{
\dfrac{(B'-B)\sqrt{\lambda_{\min}(K_T)}}{\sigma_0}
- \sqrt{T}
}^2
}.
\end{aligned}
\label{eq: final_MSE_FTP bound}
\end{equation}
When $\sigma_0=0$ and $B'\ge B$, a feasible zero-loss fit exists, so every selected minimizer satisfies $\hat f_i^{(c)}(t_j)=f_i(t_j)$. Consequently, only privacy noise contributes to the grid MSE, giving $\operatorname{MSE}_{\mathrm{FTP}}=4(B')^2\operatorname{tr}(K_T)/(n^2\mu^2T)$.
\end{proof}
\begin{remark}\textup{\textbf{(Minimum-Norm Interpolator Problem)}}
    Although both the FTP and PTF procedures involve the constrained optimization
    problems in Equations~\eqref{eq: F_operator} and~\eqref{eq: weight least square},
    the role of the constraint is different in the noiseless setting. For the
    individual trajectories, when $\sigma_0=0$ we have
    $x_{ij}=f_i(t_j)$. Hence, the empirical least-squares loss can attain its
    global minimum value zero. In particular, the constrained problem
    \[
    \hat f_i
    =
    \argmin_{\substack{f\in\mathcal H_K\\
    \|f\|_{\mathcal H_K}\le B'}}
    \frac{1}{T}\sum_{j=1}^T
    \left[f(t_j)-x_{ij}\right]^2
    \]
    admits $f_i$ as a feasible zero-loss solution. Consequently, the constraint
    does not affect the optimal value of the individual fitting problem.\vspace{1em}\\
    A convenient representative of the zero-loss solutions is the
    minimum-norm interpolant
    \begin{equation}
    \hat f_i^{\mathrm{MN}}
    =
    \argmin_{f\in\mathcal H_K}
    \|f\|_{\mathcal H_K}
    \quad
    \text{subject to}
    \quad
    f(t_j)=x_{ij},\qquad j=1,\ldots,T.
    \label{eq: minimum norm interpolator}
    \end{equation}
    Let $K_T=\left[K(t_j,t_l)\right]_{j,l=1}^T$ be the kernel Gram matrix, which is positive definite for distinct sampling locations under the Gaussian kernel. A direct Hilbert-space projection argument shows that the minimum-norm interpolation belongs to the finite-dimensional span of the kernel sections:
    \[
    \hat f_i^{\mathrm{MN}}(t)
    =
    k_T(t)^\top K_T^{-1}x_i,
    \qquad
    k_T(t)
    =
    \left(K(t,t_1),\ldots,K(t,t_T)\right)^\top.
    \]
    Moreover, $\|\hat f_i^{\mathrm{MN}}\|_{\mathcal H_K}^2=x_i^\top K_T^{-1}x_i$. Since $f_i$ itself is a feasible interpolation and
    $\|f_i\|_{\mathcal H_K}\le B$, the minimum-norm property gives
    $\|\hat f_i^{\mathrm{MN}}\|_{\mathcal H_K}
    \le
    \|f_i\|_{\mathcal H_K}
    \le B$. Thus, $\hat f_i^{\mathrm{MN}}$ is also feasible for the constrained
    least-squares problem and achieves its minimum loss, namely zero.
    Consequently, in the noiseless setting, we may take
    $\mathcal F_{\mathrm{FTP}}(x_i)
    =
    \hat f_i^{\mathrm{MN}}$. This observation allows the individual fitting step in the FTP procedure to be represented in closed form through the kernel Gram matrix.
\end{remark}
\begin{proof}\textbf{of Theorem~\ref{thm: PTF MSE lower bound}.}
For each $i$, we can still define the minimum norm interpolant as
\[
g_i
=
\argmin_{g\in\mathcal H_K}
\left\{
\|g\|_{\mathcal H_K}:
g(t_j)=f_i(t_j),\ j=1,\ldots,T
\right\},
\]
and its transformed coefficients as in the proof for Theorem~\ref{thm: FTP MSE RKHS}:
\[
\mathbf \beta_i
:=
K_T^{-1/2}\mathbf f_i,\quad
\|g_i\|_{\mathcal H_K} = \|\beta_i\|_2\le B.
\]
Denote $a_j = \sqrt{K(t_j,t_j)}$, by the RKHS reproducing property, 
\[
|f_i(t_j)|
\leq
\|f_i\|_{\mathcal{H}_K}\sqrt{K(t_j,t_j)}
\leq
Ba_j
\leq 
B' a_j.
\]
Since we are considering the noiseless case $\sigma_0=0$ and $x_i = \mathbf{f}_i$, with
\[
c_j = B'a_j.
\]
the clipping operation $\operatorname{clip}(x_{ij},c_j)$ is inactive and
\[
\bar x^{(\mathrm{clip})}
=
\bar{x}
=
\dfrac{1}{n}\sum_{i=1}^n x_i
=
\bar{\mathbf{f}}.
\]
It follows that the released private vector can be written as
\[
\widetilde x
=
\bar{\mathbf{f}}+Z
=
K_T^{1/2}\bar{\beta} + Z
,
\]
where
\[
Z\sim N(0,\tau^2W^{-1}),
\quad
\tau^2=\frac{4B'^2T}{n^2\mu^2},
\quad
W = \mathrm{diag}\{a_1^{-2},\ldots,a_T^{-2}\}.
\]
By the representer theorem, the constrained estimator has the form
\[
\hat f_{\mathrm{PTF}}^{(c)}(t)
=
\sum_{j=1}^T
(\hat\alpha^{(c)})_jK(t_j,t),
\]
\[
\hat\alpha_{\mathrm{PTF}}^{(c)}
=
\argmin_{\alpha^\top K_T\alpha\le {B'}^2}
\dfrac{1}{T} 
(\widetilde{x} - K_T\alpha)^\top W (\widetilde{x} - K_T\alpha).
\]
Denote $\beta = K_T^{1/2}\alpha$, we can rewrite the constrained problem as
\[
\hat\beta_{\mathrm{PTF}}^{(c)}
=
\argmin_{\|\beta\|_2\le B'}
(\widetilde{x} - K_T^{1/2}\beta)^\top W (\widetilde{x} - K_T^{1/2}\beta),
\quad 
\hat\beta_{\mathrm{PTF}}^{(c)} = K_T^{1/2}\hat\alpha_{\mathrm{PTF}}^{(c)}.
\]
We then analyze the corresponding unconstrained WLS estimator.\\
Define 
\[
\Phi:=K_T^{1/2},
\quad
M:=\Phi^\top\Phi = K_T,
\quad
A:= \Phi^\top W\Phi
=K_T^{1/2}WK_T^{1/2}.
\]
By
definition, the unconstrained problem is:
\[
\hat{\beta}_{\mathrm{PTF}}^{(u)}
=
\argmin_{\beta\in \R ^T} (\widetilde{x} - K_T^{1/2}\beta)^\top W (\widetilde{x} - K_T^{1/2}\beta).
\]
with the solution
\[
\widetilde\beta_{\mathrm{PTF}}^{(u)}
=
A^{-1}\Phi^\top W\widetilde x
=
\Phi^{-1}\widetilde{x}.
\]
Substituting 
$\widetilde x=K_T^{1/2}\bar{\beta} + Z= \Phi\bar\beta + Z$ gives $\widetilde\beta_{\mathrm{PTF}}^{(u)}=\bar\beta+\xi$, where $\xi:=\Phi^{-1}Z$.\\
Since $Z\sim \mathcal{N}(0,\tau^2W^{-1})$, we obtain
\[
\xi
\sim
\mathcal{N} (0,\tau^2 \Phi^{-1}W^{-1}\Phi^{-1})
\overset{d}{=}
\mathcal{N}(0,\tau^2 A^{-1}).
\]
Hence
\begin{equation}
\mathbb E
\left\|
\xi
\right\|_2^2
=
\tau^2\operatorname{tr}(A^{-1}),\quad
\mathbb E
\left\|
Z
\right\|_2^2
=
\tau^2\operatorname{tr}(W^{-1}).
\label{eq: PTF unconstrained MSE_RKHS}
\end{equation}
We next account for the constraint $\|\beta\|_2\le B'$. 
Let $\mathcal A=\{\|\widetilde\beta_{\mathrm{PTF}}^{(u)}\|_2>B'\}$ denote the event that the unconstrained estimator violates the constraint. On the complementary event $\mathcal A^c$, the unconstrained estimator is feasible and hence
$\widetilde\beta_{\mathrm{PTF}}^{(c)}
=
\widetilde\beta_{\mathrm{PTF}}^{(u)}$. Therefore
\[
\begin{aligned}
\mathbb E\|\Phi(\widetilde\beta_{\mathrm{PTF}}^{(c)}-\bar\beta)\|_2^2
&=
\mathbb E[
\|\Phi(\widetilde\beta_{\mathrm{PTF}}^{(u)}-\bar\beta)\|_2^2
\mathbf 1_{\mathcal A^c}]
+
\mathbb E[
\|\Phi(\widetilde\beta_{\mathrm{PTF}}^{(c)}-\bar\beta)\|_2^2
\mathbf 1_{\mathcal A}]
\\
&\geq
\mathbb E[
\|\Phi\xi\|_2^2
\mathbf 1_{\mathcal A^c}]
\\
&=
\mathbb E\|Z\|_2^2
-
\mathbb E[
\|Z\|_2^2
\mathbf 1_{\mathcal A}].
\end{aligned}
\]
Thus, by Equation~\eqref{eq: PTF unconstrained MSE_RKHS},
\begin{equation}
\begin{aligned}
    \operatorname{MSE}_{\mathrm{PTF}}
    =
    \dfrac{1}{T}\mathbb{E} 
    \|\widetilde{\mathbf{f}}_{\mathrm{PTF}} - \bar{\mathbf{f}}\|_2^2
    =
    \dfrac{1}{T}\mathbb{E} 
    \|\Phi(\widetilde\beta_{\mathrm{PTF}}^{(c)}-\bar\beta)\|_2^2
    \geq
    \dfrac{1}{T}\SqBra{\tau^2\operatorname{tr}(W^{-1})
    -
    R_{\mathrm{PTF}}},
\end{aligned}
\label{eq: PTF remainder decomposition_RKHS}
\end{equation}
where
\[
R_{\mathrm{PTF}}
=
\mathbb E[
\|Z\|_2^2
\mathbf 1_{\mathcal A}].
\]
It remains to bound the remainder term. Since
\[
\|\bar\beta\|_2\le \dfrac{1}{n}\sum_{i=1}^n\|\beta_i\|_2\leq B,
\]
the event $\mathcal A$ implies
\[
\|\xi\|_2
=
\|\widetilde\beta_{\mathrm{PTF}}^{(u)}-\bar\beta\|_2
\geq
\|\widetilde\beta_{\mathrm{PTF}}^{(u)}\|_2-\|\bar\beta\|_2
>
B'-B.
\]
Consequently,
\[
\mathcal A
\subseteq
\left\{
\|\xi\|_2>B'-B
\right\}.
\]
Write $\xi=\tau A^{-1/2}z,\ z\sim N(0,I_T)$. Since
\[
\|A^{-1/2}z\|_2^2
\leq
\lambda_{\max}(A^{-1})\|z\|_2^2
=
\frac{\|z\|_2^2}
{\lambda_{\min}(A)},
\]
we obtain
\[
\Pr(\mathcal A)
\leq
\Pr\left(
\|z\|_2^2
>
\frac{(B'-B)^2\lambda_{\min}(A)}
{\tau^2}
\right).
\]
Because $\|z\|_2^2\sim\chi_T^2$, apply the chi-square tail bound in Equation~\eqref{eq: chi_square bound}, we get
\begin{equation}
\Pr(\mathcal A)
\leq
\exp\left[
-\frac12
\left\{
\frac{(B'-B)\sqrt{\lambda_{\min}(A)}}{\tau}
-\sqrt T
\right\}^2
\right].
\label{eq: PTF active probability_RKHS}
\end{equation}
Finally, by the Cauchy--Schwarz inequality,
\[
R_{\mathrm{PTF}}
\leq
\left(
\mathbb E\|Z\|_2^4
\right)^{1/2}
\Pr(\mathcal A)^{1/2}.
\]
For a centered Gaussian vector with covariance
$\tau^2W^{-1}$,
\[
\mathbb E\|Z\|_2^4
=
\tau^4
\left[
\{\operatorname{tr}(W^{-1})\}^2
+
2\operatorname{tr}(W^{-2})
\right]
\le 
\tau^4 [3\{\operatorname{tr}(W^{-1})\}^2].
\]
Therefore
\[
\left(
\mathbb E\|Z\|_2^4
\right)^{1/2}
\leq
\tau^2
\sqrt{3}\operatorname{tr}(W^{-1}).
\]
Combining this inequality with Equation~\eqref{eq: PTF active probability_RKHS}
and Equation~\eqref{eq: PTF remainder decomposition_RKHS}, we prove 
Equation~\eqref{eq: PTF MSE lower bound}.
\end{proof}

\begin{proof}
\textbf{of Corollary~\ref{col: MSE ratio in RKHS}.}
Write $\tau_n = \tau =2B_n'\sqrt T/(n\mu)$. Since $B_n'$ is a non-decreasing sequence  satisfying $B_n'>B$ for sufficiently large $n$, we have
\[
\frac{B_n'-B}{\tau_n}
=
\frac{n\mu}{2\sqrt T}
\left(1-\frac{B}{B_n'}\right)
\longrightarrow\infty.
\]
Thus, the radius condition in Theorem~\ref{thm: PTF MSE lower bound} holds for all
sufficiently large $n$, and its exponential remainder tends to zero. Combining that lower bound with the exact noiseless FTP error in Theorem~\ref{thm: FTP MSE RKHS} gives
\[
\limsup_{n\to\infty}
\frac{\operatorname{MSE}(\widetilde f_{\mathrm{FTP}})}
     {\operatorname{MSE}(\widetilde f_{\mathrm{PTF}})}
\le
\frac{\operatorname{tr}(K_T)}
     {T\operatorname{tr}(W^{-1})}.
\]
By the definition
$W=\operatorname{diag}\{K(t_j,t_j)^{-1}:j=1,\ldots,T\}$,
\[
\operatorname{tr}(W^{-1})
=\sum_{j=1}^T K(t_j,t_j)
=\operatorname{tr}(K_T).
\]
Hence, we get the claimed bound of $1/T$.
\end{proof}
\section{Additional Simulation Results}
\label{sec: Additional Simulation Results}
In this section, we provide the simulation setup and results for polynomial spaces and supplementary numerical results for the RKHS setting. 

\textbf{Polynomial Setup.} Consider the functional model: $f_1,\ldots,f_n\in \mathcal{H}_m$ with $L^2[0,1]$ norm and inner product, where $\mathcal{H}_m=\Set{\theta^\top \phi(t):\theta\in \R^m, t\in [0,1]}$, and the orthonormal polynomial basis $\phi(t) = (\phi_1(t),\ldots,\phi_m(t))^\top$ satisfies that the polynomial degree of $\phi_k(t)$ is $k-1$ for $k=1,\ldots,m$. For example, when $m=2$, the orthonormal basis is $(\phi_1(t),\phi_2(t))^\top = (1,\sqrt{12}(t-0.5))^\top$. We construct this basis iteratively based on shifted orthonormal Legendre basis. Then, for each $i$, we generate $f_i(t) = \theta_i^\top \phi(t)$ with $\theta_i$ i.i.d. in the $m$-dimensional ball with radius $B=1$, which ensures $\|f_i\|_{L^2} = \|\theta_i\|_2\leq B$. The data matrix $X$ is generated based on the model: $x_{ij} = \theta_i^\top \phi(t_j) + \epsilon_{ij}, \epsilon_{ij}\overset{i.i.d.}{\sim} N(0,\sigma_0^2)$, with an equally spaced location grid between 0 and 1. Specifically, $t_j = (j-1)/(T-1), j=1,\ldots,T$ and we fix $T=20$. Then we obtain the two private estimators $\widetilde{f}_\mathrm{FTP}$ and $\widetilde{f}_\mathrm{PTF}$ through Algorithm~\ref{alg: FTP} and Algorithm~\ref{alg: PTF}, respectively. Finally, we use the IMSE metric to evaluate the two estimators compared to the sample mean function $\bar{f} = n^{-1}\sum_{i=1}^n f_i$. We try different settings for sample sizes, basis dimensions, privacy budgets, constraint radius, and observation noise standard deviations, where $n\in \Set{2^2, \ldots, 2^{12}}$, $m\in \Set{2,4,6,8}$, $\mu\in \Set{0.1, 1}$, $B'\in\Set{B, 1.1B, \sqrt{\log n}, \log n, \sqrt{n}, n^{2/3}}$, and $\sigma_0\in\Set{0, 0.1, 1, 10}$, respectively. 

\begin{figure}[t]
    \centering
    \includegraphics[width=\linewidth]{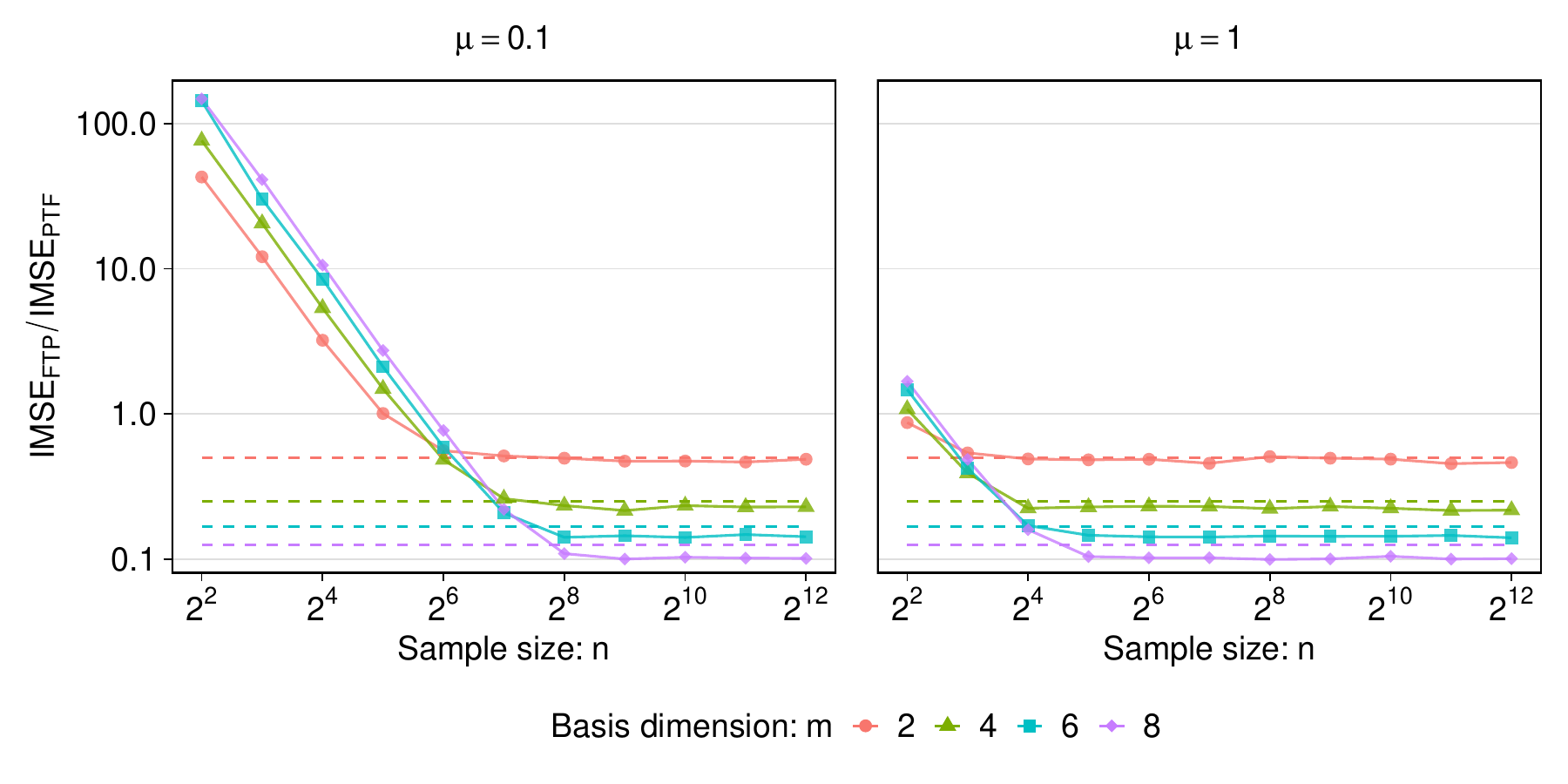}
    \vspace{-6mm}
    \caption{Empirical IMSE ratio of FTP to PTF versus sample size $n$ in polynomial spaces, with $\sigma_0=0$, $B'=1.1$. The panels correspond to $\mu=0.1$ and $\mu=1$. Colors and symbols identify the basis dimension $m$. Dashed lines show the corresponding $1/m$ asymptotic upper bounds. Both axes use logarithmic scales.}
    \label{fig: IMSE Ratio with Bound in Polynomial}
\end{figure}

Figure~\ref{fig: IMSE Ratio with Bound in Polynomial} examines the noiseless IMSE comparison between FTP and PTF in the polynomial space with $B'=1.1$. The IMSE ratios generally decrease with sample size and stabilize near or below the corresponding $1/m$ reference levels, consistent with Corollary~\ref{col: IMSE Ratio for Polynomials}. Larger basis dimensions yield smaller ratios at large sample sizes, while stronger privacy protection ($\mu=0.1$) delays the transition to this regime. At small sample sizes, PTF can outperform FTP, particularly under stronger privacy protection.

\begin{figure}[htbp]
    \centering
    \includegraphics[width=\linewidth]{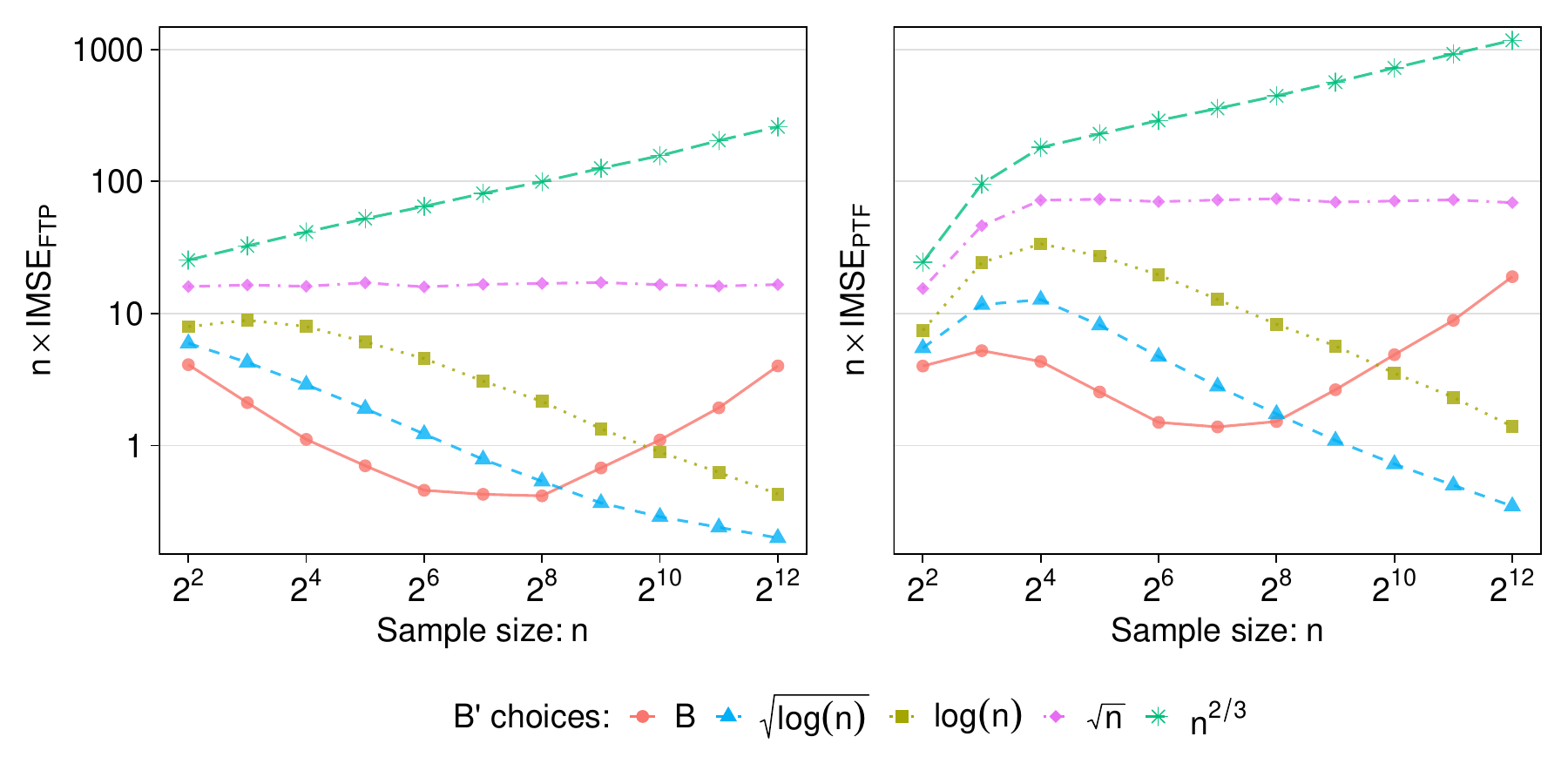}
    \vspace{-2em}
    \caption{Scaled empirical IMSE ($n\times\operatorname{IMSE}$) for FTP (left) and PTF (right) in polynomial spaces under different fitting radii $B'$, with $\sigma_0=1$, $\mu=1$, $m=4$. Both axes use logarithmic scales.}
    \label{fig: Poly_Bprime_choices}
\end{figure}

\begin{figure}[htbp]
    \centering
    \includegraphics[width=\linewidth]{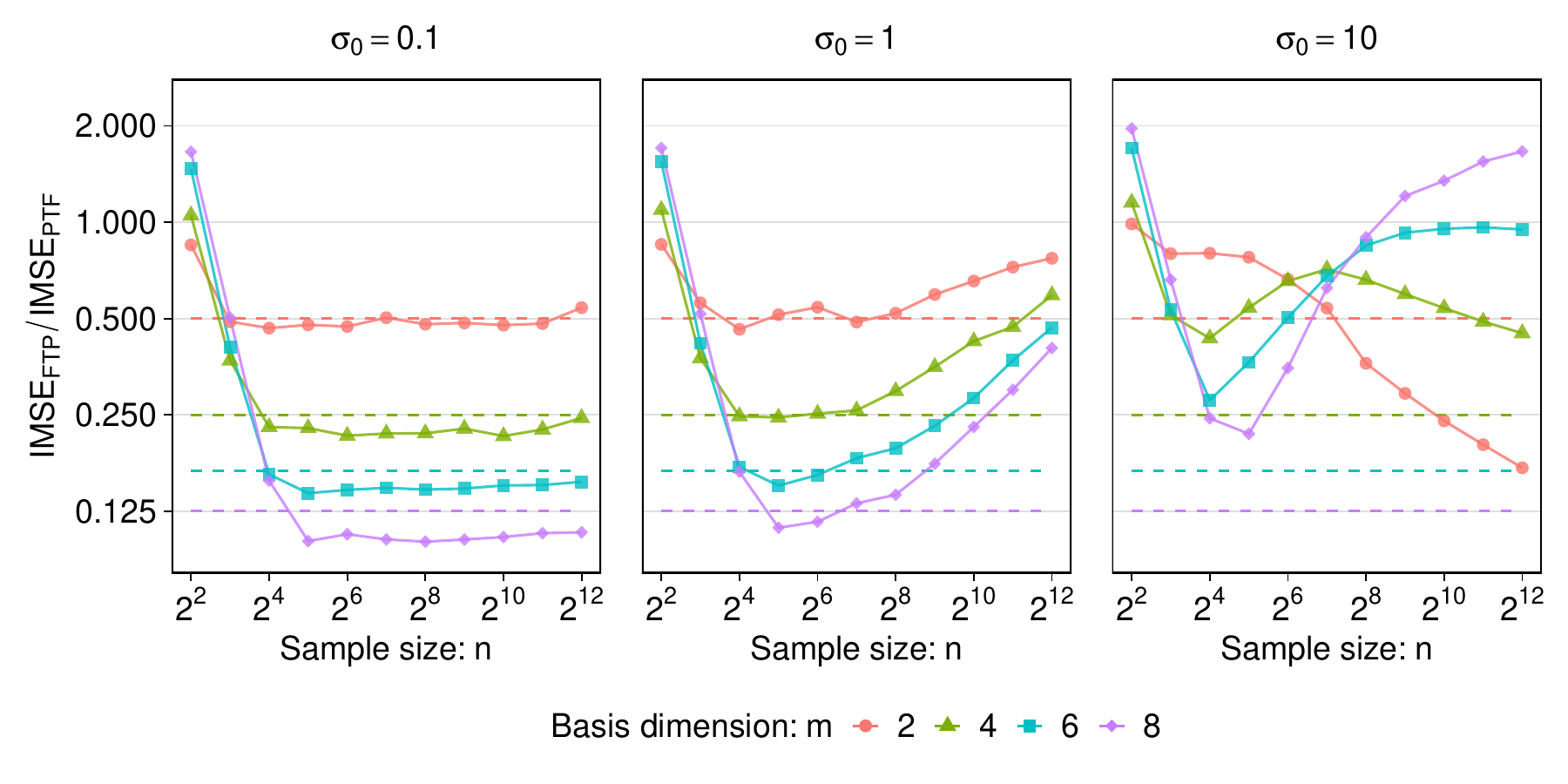}
    \vspace{-2em}
    \caption{Empirical IMSE ratio of FTP to PTF versus sample size $n$ in polynomial spaces for $\sigma_0\in\{0.1,1,10\}$, with $\mu=1$, $B'=\sqrt{\log n}$. Colors and symbols identify the basis dimension $m$. Dashed lines show the noiseless asymptotic reference bounds $1/m$. Both axes use logarithmic scales.}
    \label{fig: IMSE Ratio with Bound in Polynomial with noise}
\end{figure}

Figure~\ref{fig: Poly_Bprime_choices} compares the scaled IMSE under different fitting radii with $\sigma_0=1$, $\mu=1$, and $m=4$. The oracle radius $B'=B$ performs well at intermediate sample sizes, but its scaled error increases at larger $n$. Among the displayed choices, $B'=\sqrt{\log n}$ gives the smallest error for both pipelines at the largest sample sizes. The approximately flat curves for $B'=\sqrt n$ and increasing curves for $B'=n^{2/3}$ illustrate the privacy cost of allowing the fitting radius to grow rapidly. These results complement the RKHS experiments by showing that the preferred radius over the simulated range can depend on the fitting space.

Figure~\ref{fig: IMSE Ratio with Bound in Polynomial with noise} compares the IMSE ratios under different observation-noise levels with $\mu=1$ and $B'=\sqrt{\log n}$. When $\sigma_0=0.1$, the ratios remain close to or below the noiseless reference levels at larger sample sizes. When $\sigma_0=1$, the ratios initially decrease and then increase, although FTP retains an advantage at the largest sample size for every displayed basis dimension. For $\sigma_0=10$, the behavior depends more strongly on $m$: PTF becomes preferable for $m=8$ at larger sample sizes, whereas FTP remains preferable for the smaller dimensions at the end of the displayed range. Thus, the polynomial experiments support FTP's advantage under low observation noise while showing that its robustness to larger noise levels depends on the fitting dimension.

\begin{table}[t]
\centering
\caption{Noiseless RKHS simulation results for $B=1$, $B'=1.1$, and Gaussian-kernel bandwidth $\rho=0.2$. The ratio is the Monte Carlo MSE of FTP divided by that of PTF. $\widehat p_{\mathrm{PTF}}$ is the empirical probability that the PTF fitting constraint is active over 1000 repetitions; a displayed zero does not imply a theoretically zero probability.}
\label{tab: RKHS noiseless details}
\small
\setlength{\tabcolsep}{3pt}
\begin{tabular}{rr rrrr rrrr}
\toprule
$T$ & $n$ & \multicolumn{2}{c}{$\mathrm{MSE}_{\mathrm{FTP}}$} & \multicolumn{2}{c}{$\mathrm{MSE}_{\mathrm{PTF}}$} & \multicolumn{2}{c}{Ratio} & \multicolumn{2}{c}{$\widehat p_{\mathrm{PTF}}$} \\
\cmidrule(lr){3-4}\cmidrule(lr){5-6}\cmidrule(lr){7-8}\cmidrule(lr){9-10}
& & $\mu=0.1$ & $\mu=1$ & $\mu=0.1$ & $\mu=1$ & $\mu=0.1$ & $\mu=1$ & $\mu=0.1$ & $\mu=1$ \\
\midrule
2 & $2^{2}$ & $2.94\times 10^{1\phantom{-}}$ & \num{3.15e-01} & \num{9.05e-01} & \num{3.40e-01} & 32.427 & 0.927 & 0.991 & 0.575 \\
 & $2^{4}$ & $1.85\times 10^{0\phantom{-}}$ & \num{1.97e-02} & \num{6.83e-01} & \num{3.38e-02} & 2.710 & 0.582 & 0.868 & 0.090 \\
 & $2^{6}$ & \num{1.11e-01} & \num{1.19e-03} & \num{1.58e-01} & \num{2.33e-03} & 0.704 & 0.509 & 0.407 & 0.000 \\
 & $2^{8}$ & \num{7.31e-03} & \num{7.47e-05} & \num{1.40e-02} & \num{1.49e-04} & 0.523 & 0.502 & 0.029 & 0.000 \\
 & $2^{10}$ & \num{4.83e-04} & \num{4.67e-06} & \num{8.82e-04} & \num{9.05e-06} & 0.547 & 0.517 & 0.000 & 0.000 \\
 & $2^{12}$ & \num{2.81e-05} & \num{2.77e-07} & \num{5.66e-05} & \num{5.69e-07} & 0.496 & 0.486 & 0.000 & 0.000 \\
\midrule
4 & $2^{2}$ & $3.04\times 10^{1\phantom{-}}$ & \num{3.02e-01} & \num{5.38e-01} & \num{3.16e-01} & 56.554 & 0.954 & 1.000 & 0.946 \\
 & $2^{4}$ & $1.86\times 10^{0\phantom{-}}$ & \num{1.98e-02} & \num{4.50e-01} & \num{6.05e-02} & 4.140 & 0.328 & 0.998 & 0.379 \\
 & $2^{6}$ & \num{1.14e-01} & \num{1.18e-03} & \num{1.97e-01} & \num{4.64e-03} & 0.576 & 0.255 & 0.811 & 0.011 \\
 & $2^{8}$ & \num{7.53e-03} & \num{7.38e-05} & \num{2.59e-02} & \num{2.93e-04} & 0.291 & 0.251 & 0.235 & 0.000 \\
 & $2^{10}$ & \num{4.70e-04} & \num{4.69e-06} & \num{1.81e-03} & \num{1.79e-05} & 0.259 & 0.262 & 0.000 & 0.000 \\
 & $2^{12}$ & \num{2.78e-05} & \num{2.88e-07} & \num{1.14e-04} & \num{1.17e-06} & 0.244 & 0.247 & 0.000 & 0.000 \\
\midrule
6 & $2^{2}$ & $2.90\times 10^{1\phantom{-}}$ & \num{2.99e-01} & \num{5.24e-01} & \num{2.85e-01} & 55.296 & 1.048 & 1.000 & 1.000 \\
 & $2^{4}$ & $1.89\times 10^{0\phantom{-}}$ & \num{1.90e-02} & \num{4.35e-01} & \num{6.09e-02} & 4.343 & 0.312 & 1.000 & 0.891 \\
 & $2^{6}$ & \num{1.22e-01} & \num{1.19e-03} & \num{1.93e-01} & \num{6.62e-03} & 0.630 & 0.180 & 0.996 & 0.141 \\
 & $2^{8}$ & \num{7.36e-03} & \num{7.38e-05} & \num{2.96e-02} & \num{4.47e-04} & 0.249 & 0.165 & 0.741 & 0.000 \\
 & $2^{10}$ & \num{4.57e-04} & \num{4.72e-06} & \num{2.79e-03} & \num{2.82e-05} & 0.164 & 0.167 & 0.004 & 0.000 \\
 & $2^{12}$ & \num{2.81e-05} & \num{2.85e-07} & \num{1.73e-04} & \num{1.74e-06} & 0.162 & 0.164 & 0.000 & 0.000 \\
\midrule
8 & $2^{2}$ & $3.03\times 10^{1\phantom{-}}$ & \num{3.06e-01} & \num{5.16e-01} & \num{2.73e-01} & 58.788 & 1.121 & 1.000 & 1.000 \\
 & $2^{4}$ & $1.95\times 10^{0\phantom{-}}$ & \num{1.88e-02} & \num{4.29e-01} & \num{5.61e-02} & 4.537 & 0.335 & 1.000 & 1.000 \\
 & $2^{6}$ & \num{1.20e-01} & \num{1.19e-03} & \num{1.83e-01} & \num{6.61e-03} & 0.656 & 0.180 & 1.000 & 0.918 \\
 & $2^{8}$ & \num{7.38e-03} & \num{7.06e-05} & \num{2.73e-02} & \num{5.68e-04} & 0.270 & 0.124 & 0.995 & 0.192 \\
 & $2^{10}$ & \num{4.56e-04} & \num{4.78e-06} & \num{2.98e-03} & \num{3.70e-05} & 0.153 & 0.129 & 0.755 & 0.000 \\
 & $2^{12}$ & \num{2.88e-05} & \num{2.80e-07} & \num{2.26e-04} & \num{2.32e-06} & 0.127 & 0.121 & 0.024 & 0.000 \\
\bottomrule
\end{tabular}
\end{table}

Table~\ref{tab: RKHS noiseless details} supplements Figure~\ref{fig: MSE Ratio with Bound in RKHS} with the individual MSEs and the empirical probability that the PTF fitting constraint is active. This probability generally decreases with sample size, but remains higher for larger $T$ or smaller $\mu$. For example, with $T=8$ and $n=2^{10}$, the estimated probability is $0.755$ for $\mu=0.1$, compared with zero observed activations for $\mu=1$. As constraint activation becomes less frequent, the error ratios approach the corresponding $1/T$ reference levels, up to Monte Carlo variation.

\vskip 0.2in

\clearpage
\bibliography{references}

\end{document}